\documentclass{article}
\usepackage[toc,page]{appendix}
\usepackage{verbatim}
\usepackage{amssymb}
\usepackage{bm}
\usepackage{graphicx}
\usepackage{amsmath}
\DeclareMathAlphabet{\mathpzc}{OT1}{pzc}{m}{it}
\newcommand{\sG}{\mathpzc{G}}
\newcommand{\sT}{\mathpzc{T}}
\newtheorem{theorem}{Theorem}

\newtheorem{algorithm}[theorem]{Algorithm}

\newtheorem{definition}[theorem]{Definition}
\newtheorem{example}[theorem]{Example}

\newtheorem{lemma}[theorem]{Lemma}

\newtheorem{proposition}[theorem]{Proposition}
\newtheorem{procedure}[theorem]{Procedure}
\newtheorem{remark}[theorem]{Remark}

\newenvironment{proof}[1][Proof]{\textbf{#1.} }{\ \rule{0.5em}{0.5em}}
\newenvironment{pf}[1][Proof of Theorem]{\textbf{#1} }{\ \rule{0.5em}{0.5em}}

\newcommand{\supp}{\mathrm{supp}}

\newcommand{\Sym}{\mathrm{Sym}}

\newcommand {\tr}{\mathrm{tr}}

\newcommand{\so}{\mathfrak{so}}

\newcommand{\h}{\mathfrak{h}}

\newcommand{\gl}{\mathfrak{gl}}

\renewcommand{\O}{\mathcal{O}}

\newcommand{\F}{\mathcal{F}}
\newcommand{\SO}{\mathrm{SO}}

\renewcommand{\P}{\mathbf{P}}
\newcommand{\U}{\mathrm{U}}

\newcommand{\SL}{\mathrm{SL}}
\renewcommand{\sl}{\mathfrak{sl}}

\newcommand{\sd}{\partial}

\newcommand{\ity}{_{\infty}}
\newcommand{\g}{\mathfrak{g}}
\renewcommand{\b}{\mathfrak{b}}
\newcommand{\n}{\mathfrak{n}}
\newcommand{\ddef}{\mathrm{def}}
\newcommand{\Maps}{\mathrm{Maps}}
\newcommand{\End}{\mathrm{End}}
\newcommand{\Hom}{\mathrm{Hom}}

\begin{document}
\title{On quasimaps to quadrics} 

\author{}
\date{\today}
\maketitle
\begin{abstract}
We put on a rigorous basis results of Aisaka and  Aldo Arroyo on curved $\beta\gamma$ systems on a quadric. Drinfeld's quasimaps turned to be indispensable in accurate computation of mini-BRST cohomology. 
\end{abstract}
\section{Introduction}\label{S:introduction}
A mathematically  rigorous derivation of the Hilbert space of string theory in the formalism of pure spinors \cite{Berkovits} remains a challenging problem \cite{AABN}. The difficulty is to give a mathematical meaning  a $\beta\gamma$-system with the target  the space of pure spinors $\SO(10)/\U(5)$, more precisely an affine cone over it.The direct attack on the problem  would  require to develop  analysis on the space of loops  in  the target. The paper \cite{AA} proposes to investigate   a more simple model where the target is an affine quadric. 
 The space of loops in the quadric is an infinite-dimensional "manifold", defined by an infinite set of equations. The Koszul or mini-BRST(in terminology of \cite{AA}) complex, constructed from the set of equations hide all difficulties of analysis in the BRST differential. The advantage of this formulation is that constraint fields of the original problem are replaced by free fields  of mini-BRST setup.

More precisely  following \cite{AA} we consider a space of loops that lie on a quadric in a complex linear space $V$. It is convenient to choose an orthonormal  basis $e_i,i=1\dots,n$. Then a loop \[\lambda:S^1\rightarrow V\] is completely characterizes  by $n$ functions $\lambda^i(z),z\in S^1\subset \mathbb{C}$: \[\lambda(z)=\lambda^1(z)e_1+\cdots+\lambda^n(z)e_n\]
The loop belongs to the quadric if $\sum_{i=1}^n(\lambda^i(z))^2=0$ for all $z$. One can go on and expand $\lambda$ into Fourier  series \[\lambda(z)=\sum_{k\in \mathbb{Z}}\sum_{i=1}^n\lambda^i[k]e_iz^k\]
From this we obtain an infinite set of equations on Fourier coefficients:
\begin{equation}\label{E:reltrivial}
f[k]=\sum_{t+s=k}\sum_{i=1}^n\lambda^i[t]\lambda^i[s]=0
\end{equation}
Following the standard in the gauge theory  and the commutative algebra construction  we define Koszul complex. For every constraint $f_k$ we  introduce an odd variable $c_k$, that carries degree one. In the space of functions in $\lambda^k_i,c_k$ we define a differential by the formula \[\sum_kf_k\frac{\partial}{\partial c_k}.\]
One of the assertions that has been used implicitly many times in \cite{AA} was that the complex has  cohomology only in degree zero. Though the authors attribute this to irreducibility of the constraint, this statement hard to make precise in infinite dimensions because cohomology depend also on the class of functions we are dealing with.

One of the goals of the present work is to make this statement precise.

The key technical novelty is that we use  finite-dimensional spaces  of Drinfeld's quasimaps to approximate the loop spaces.\footnote{As author learned from \cite{Braverman} this idea was first introduced by Drinfeld in a context of quasimaps to full flag spaces. } This means that we  work with 
Fourier  series of the form  \[\lambda(z)=\sum_{-N_1\leq k\leq N_2}\sum_{i=1}^n\lambda^i[k]e_iz^k\]
Equations $f_k=0$ define a finite-dimensional algebraic variety, which typically is singular. As a side remark we mention that De Rham and Dolbeault  complexes have a limited usefulness in this setup. Never the less we can study algebras of polynomial  functions $AQ=AQ_{-N_1}^{N_2}(n)$ in  variables $\lambda^i[k]$ subject to relations $f[k]$ and define Koszul complexes $C=C_{-N_1}^{N_2}(n)$ using even $\lambda^i[k]$ $N_1\leq k\leq N_2$ and odd $c[k]$ $-2N_1\leq k\leq 2N_2$ as generators. One of our results is
\begin{theorem}\label{T:main}
The zero cohomology of the  complex $C_{-N_1}^{N_2}(n)$ is equal to  $AQ_{-N_1}^{N_2}(n)$, $n\geq 3$. All other cohomology vanish.
\end{theorem}
In our paper we use  ideas of  \cite{Hodge}
 in the version of   \cite{SottileSturmfels} \cite{Ruffo}  to compute Poincar\'{e} series of $AQ$.  The method lets us not only to compute the series but to prove important Koszul property of algebra $AQ$. This becomes crucial in the proof acyclicity of  resolution proposed by Aisaka and Arroyo.

We believe  that direct limit of $Spec(AQ_{-N_1}^{N_2})$  capture most of the pertinent properties of loop spaces for "good"
 targets. This lets us to avoid difficulties of infinite-dimensional analysis and allows to use methods of commutative algebra and algebraic geometry.   We want to emphasize  that the method can easily be  extended to  a more realistic case of pure spinors. We plan to elaborate on this in our following publication.

Few questions were left unanswered:
\begin{itemize}
\item Do spaces of Drinfeld's quasimaps have intrinsic physical meaning or it is just a convenient mathematical tool?
\item It is known (see e.g. \cite{Nekrasov})  that $\beta\gamma$-systems might suffer from anomalies.  How anomalies appear in the proposed approach? (See a remark in the end of Section \ref{S:quasimaps})
\item Depending on the answer of the previous item, what is the construction of the Virasoro action?
\item Characterize in geometric terms conical algebraic varieties whose $\beta\gamma$-systems systems satisfy $*$-duality. 
Give a geometric construction   of $*$ duality using the language of quasimaps. 
\end{itemize}

We hope to address these questions in the near future.

The paper is organized as  follows. In Section \ref{quadric} we illustrate our methods of computation of Poincar\'{e} series with a simplest example of a quadric.
Section \ref{S:quasimaps} is technically central, where we prove straightened law and Koszul property for $AQ$. In Section \ref{S:generating} compute Poincar\'{e} series of $AQ$. In Section \ref{S:Hilbert} we compute semi-infinite  cohomology of      Lie algebra $LH$. This cohomology will serve an approximation to semi-infinite cohomology of algebra $LH^{\infty}_{-\infty}$ discussed in Section \ref{S:limits}. The latter is an algebraic counterpart of quantum mini-BRST complex of \cite{AA}.

The paper is supplemented by two appendices. In the first we discuss exceptional cases of quadrics in 4,3 or 2-D spaces. The basics of theory of Gr\"{o}bner bases is reviewed in the second appendix. 

{\bf \large Acknowledgment}

The author would like to thank N.Berkovits, M. Finkelberg, V. Gorbunov, A.Schwarz  for stimulating discussions. The author is also thankful to Max-Plank-Institute (Bonn) for excellent recearch environment, where part of this work has been done. 

\section{Finite-dimensional quadric}\label{quadric}
Before we delve into algebra of quasimaps  let us treat the simplest case of maps of degree zero to a quadric, that is   the quadric itself . The algebra of algebraic functions on an affine non-degenerate quadric $A=\bigoplus_{i\geq 0} A_i$ is a quotient $\mathbb{C}[\lambda^i]/(\sum_{i=1}^n \lambda^2_i)$  has Poincar\'{e} series  \[A(t)=\sum_{i\geq 0} \dim A_it^i=\frac{1-t^2}{(1-t)^{\dim(V)}}.\]
Let us see how our method reproduces this result.




\subsection{The case of an even-dimensional $V$}\label{S:even}
We suppose  that $\dim(V)=2n$, $n\geq 3$ and decompose  \[V=W+W^*\]  into a direct sum of complementary isotropic subspaces. The space $W^*$ is dual to $W$ with the pairing defined by the formula $<g,f>=(g,f)$.
We  choose a basis $g_1,\dots,g_n$ in $W$ and  the dual basis $f^1,\dots,f^n$ in $W^*$.

The adjoint representation of  the complex Lie algebra  $\so_{2n}$ is isomorphic to \[\Lambda^2(V)\cong \Lambda^2(W^*)+W^*\otimes W +\Lambda^2(W)\]
We identify the space $W^*\otimes W $ with the Lie subalgebra $\gl_n$ of $\so_{2n}$. A subalgebra of $\gl_n$ of diagonal matrices (with respect to the basis  $g_1,\dots,g_n$)
is  a Cartan subalgebra $\h$ of $\so_{2n}$. It is spanned by $f^ig_i, i=1\dots,n$. Denote  by $\rho$ the fundamental representation  \[\rho:\so_{2n}\rightarrow \End(V)\] 
Simple positive root vectors  of $\so_{2n}$ relative to $\h$ are $r=g_1\wedge g_2$ and $r_i=g_{i+1}\otimes f_i$.
The operators $R=\rho(r),R_i=\rho(r_i)$ 
 \[
R g_k=0\quad 
Rf^k=
\begin{cases}g_2  & \text{if $k=1$,}\\
-g_1  & \text{if $k=2$,}\\
0 &\text{if $ k\neq1,2$}\
\end{cases}
\]

\[
R_{i}g_k=
\begin{cases}g_{k+1}  & \text{if $i=k$,}\\
0 &\text{if $i\neq k$}\
\end{cases}
\quad 
R_{i}f^k=
\begin{cases}-f^{k-1}  & \text{if $i=k-1$,}\\
0 &\text{if $i\neq k-1$}\
\end{cases}
\]
$i=1,\dots,n-1$
are of great importance in our construction.

%
%
We define the following Hasse diagram:
\[\sG(2n):\langle f^n \rangle \overset{R_{n-1}}\longrightarrow \dots\overset{R_2} \longrightarrow \langle f^{2} \rangle \begin{array}{c}\overset{R_1\quad \text{ }}\nearrow  \begin{array}{c} \langle f^1 \rangle \\ \text{ } \end{array} \overset{\text{ }\quad R}\searrow \\ \underset{R \quad  \text{ }}\searrow   \begin{array}{c} \text{ }\\ \langle g_1 \rangle  \end{array} \underset{\text{ }\quad R_1}\nearrow \end{array}\langle  g_2 \rangle  \overset{R_2}\longrightarrow \dots  \overset{R_{n-1}}\longrightarrow  \langle g_n \rangle \]
formed by weight subspaces of $V$.

In  this diagram the brackets $\langle \rangle$ stand for spans of vectors. One-dimensional linear spaces $\alpha,\beta$ from the above diagram are connected by an arrow $\alpha\rightarrow \beta$ if one of the operators  $R$ or $R_i$ transforms   linear space $\alpha$  isomorphically to linear space $\beta$. 

The diagram will be used  for describing a basis in the algebra \[A=\mathbb{C}[g_1,\dots,g_n, f^1,\dots,f^n]/(\sum_{i=1}^n g_if^i)\]

The graph defines a partial order on the set of generators $\sG=\sG(2n)$  of algebra $A$ ($\alpha\leq \beta$ iff there is an arrow $\alpha \rightarrow \beta$). In this  poset  any two elements $\alpha,\beta \in \sG $ have a {\it unique} supremum and infinum :  \[\delta\leq \alpha, \beta \leq \gamma\] A poset with this property is called a {\it lattice}. We denote \[\gamma\overset{\ddef}{=}\alpha\vee \beta, \quad \delta\overset{\ddef}{=}\alpha\wedge \beta.\]

It is convenient to use uniform notations for generators of $A$ \[\{e_{\alpha}|\alpha\in \sG\}=\{g_1,\dots,g_n,f^1,\dots,f^n\}\]
Thus $A$ is a quotient of $\mathbb{C}[e_{\alpha}]\quad \alpha\in \sG$\[0\rightarrow I\rightarrow \mathbb{C}[e_{\alpha}]\overset{p}\rightarrow A\rightarrow 0 \]

Any monomial $m=\prod_{i=1}^k e_{\alpha_i}\in \mathbb{C}[e_{\alpha}]$ can be rewritten as  $m=\prod_{\alpha\in B\subset \sG} e_{\alpha}^{n_{\alpha}}$. We call $B$ the support $\supp(m)$. 

A set $B\subset \sG$ is called a path if all elements of $B$ are comparable. The poset $\sG$ has precisely two incomparable  elements $\langle g_1\rangle $ and $\langle f^1\rangle$.

Observe that the  only $g_1f^1$  of the monomials of  the defining relation 
\begin{equation}\label{E:qrelation}
\sum_{i=1}^n g_if^i
\end{equation} has support that is not a path. It is useful to rewrite the defining relation in the form 
\begin{equation}\label{E:relationzero}
g_1f^1=-g_2f^2-\sum_{i\geq 3} g_if^i
\end{equation}
It is equivalent to   
\begin{equation}\label{E:relationgen}
e_{\alpha}e_{\beta}=c_{\alpha\vee \beta,\alpha\wedge \beta} e_{\alpha\vee \beta }e_{\alpha\wedge \beta }+\sum_{ \begin{array}{c}\gamma\geq \alpha\vee \beta \\\gamma' \leq \alpha \wedge \beta \end{array}}c_{\gamma,\gamma'} e_{\gamma}e_{\gamma'}
\end{equation}
with  coefficients $c_{\gamma,\gamma'}$ equal to $-1$ or $0$.


A straightforward  way to determine $A(t)$  would be to construct a basis $m^i_k$ in $A_i$, count the number of elements in this basis and write the generating function. Such  direct approach technically is not very convenient because the spaces $A_i$ are defined as quotients of graded spaces of polynomial algebra $\mathbb{C}[e_{\alpha}]$. It is much more transparent to do the counting of dimensions  directly  in terms of  $\mathbb{C}[e_{\alpha}]$. 

One can give the following restatement  of the problem.
 It is to  find a subset $E$ in the set of monomials in $\mathbb{C}[e_{\alpha}]$ such that $p(E)$ is a basis in $A$.

Not all the monomial $m\in \mathbb{C}[e_{\alpha}] $ are linearly independent in $A$ because of the relation (\ref{E:qrelation}). In particular monomial $g_1f^1$, whose support $\supp(g_1f^1)$ is not a path can be represented as a sum of monomials $-g_2f^2-\sum_{i\geq 3} g_if^i$, whose support are paths. This can be generalized to arbitrary monomials: images of  monomials whose support is a path (we call them $\sG$-monomials ) span $A$. 
\begin{proposition}
Suppose that $\dim(V)$ is $2n,n\geq 3$. Let $A$ be the algebra of algebraic  function on affine non-degenerate quadric in $V^*$. Let $\sG(2n)$ be a lattice define above. Then  the image of the set of $\sG$-monomials in $A$ is a basis.
\end{proposition}
\begin{proof}
If a linear combination of $\sG$-monomials belongs to the ideal then it is divisible by the relation (\ref{E:qrelation}). But the relation contains a non $\sG$-monomial. Thus the original linear combination contains a non $\sG$-monomial.
\end{proof}
\begin{definition}
Let $\sG$ be a lattice, $B$ be a graded algebra. Suppose that the first graded component  $B_1$ has a basis $<e_{\alpha}>$ labeled by elements of $\sG$. Suppose that in algebra $B$ relations (\ref{E:relationgen}) hold. We say that $B$ is an algebra with straightened law if the images of  $\sG$-monomials form a basis in $B$.
\end{definition}
\begin{remark}
Our definition of algebra with straightened law is a particular case of a Hodge algebra (\cite{DCEP} or \cite{Hibi}).
\end{remark}
We see that $A$ is an algebra with straightened law.

Our results  can be used to find a formula for $A(t)$. It is not hard to describe all $\sG$-monomials.
These are \[(f^n)^{\alpha_n}\cdots(f^2)^{\alpha_2} (f^1)^{\alpha_1}g_2^{\beta_2}\cdots g_n^{\beta_n}\] and \[(f^n)^{\alpha_n}\cdots(f^2)^{\alpha_2} g_1^{\beta_1}g_2^{\beta_2}\cdots g_n^{\beta_n}\]
$\alpha_i,\beta_j\geq0$

These monomials form sets $E_1$ and $E_2$.  

Let $L$ be a subset of $E$. Defines a generating function $L(t)$ as \[L(t)=\sum_{i\geq 0}\#\{m\in L|\deg(m)=i\}\]

Then \[\begin{split}&A(t)=E_1(t)+E_2(t)-(E_1\cap E_2)(t)=\\
&=\frac{1}{(1-t)^{2n-1}}+ \frac{1}{(1-t)^{2n-1}}-\frac{1}{(1-t)^{2n-2}}=\frac{1-t^2}{(1-t)^{2n}}\end{split}\]

\subsection{The case of an odd-dimensional $V$}
Let us carry out cursorily a similar  analysis of an affine non-degenerate quadric in odd-dimensional space $V$. Suppose $2n+1=\dim(V)\geq 1$. We can choose a basis $(f^1,\dots,f^n,g_1,\dots,g_n,h)$ in which  equation $q$ of the quadric is  
\begin{equation}\label{E:relationodd}
\sum_{i=1}^nf^ig_i+h^2.
\end{equation}
and \[ \langle f^1,\dots,f^n\rangle +\langle g_1,\dots,g_n\rangle +\langle h\rangle=W^*+W+\langle h\rangle\]
 The adjoint representation of $\so_{2n+1}$ is isomorphic to \[\Lambda^2(V)\cong \Lambda^2(W)+W\otimes W^*+\Lambda^2(W^*)+W\otimes \langle h\rangle+W^*\otimes \langle h\rangle\]
Positive simple root vectors   are 
$r=g_1 h$ and $r_i=g_{i+1} f_i$. Cartan subalgebra is generated by $f^ig_i$.

The corresponding diagram $\sG(2n+1)$ constructed with the aid of  $R=\rho(r),R_i=\rho(r_i)$ is 
\[\sG(2n+1):\langle f^n \rangle \overset{R_{n-1}}\longrightarrow \dots\overset{R_1} \longrightarrow \langle f^{1} \rangle \overset{R} \longrightarrow \langle h\rangle \overset{R} \longrightarrow \langle  g_1 \rangle \overset{R_1}\longrightarrow \dots  \overset{R_{n-1}}\longrightarrow \langle g_n\rangle \]
One of the distinctions  of even and odd dimensional cases is that in odd case the Weil group $W$ of $\SO(2n+1)$ does not act transitively on the set of weight vectors $f^i,g_j,h$. This is because the weight of $h$ is zero, the weights of $f^i,g_j$ are non zero. This makes it impossible for $h$ to belong to the same orbit of the Weil group because  the group acts linearly on  weights. 

A representation of semi-simple group $G$ is called {\it minuscule} if $W$ acts transitively on the set of weights. Thus the fundamental representation of $\so_k$ is minuscule only  for even $k$.

To extend a definition of algebras with straightened laws to the case of odd quadric we introduce the following definition.
\begin{definition}\label{D:wlattice}
We call  $\sT$ a weakly  partly ordered set if
it has  a (partly defined) binary relation $\leq$  which is antisymmetric, and transitive, i.e., for all $a, b,$ and $c$ in $\sT$, we have that:
 \begin{itemize}
    \item if $a \leq b$ and $b \leq a$ then $a = b$ (antisymmetry);
    \item if $a \leq b$ and $b \leq c$ then $a \leq c$ (transitivity).
\end{itemize}
\end{definition}
The above definition differs from the standard definition  of  a partly ordered set in the item of reflexivity. For a weak  partial  order condition $a\leq a$ might fail for some exceptional  $a$.
\begin{definition}
We call weakly  partly ordered set $\sT$ a weak lattice  if for any two elements $\alpha,\beta$ there there are unique supremum $\alpha\vee \beta,$ and infinum,  $\alpha\wedge \beta.$
\end{definition}

We define a structure of a weak lattice on the set of vertices of the diagram $\sG(2n+1)$ using the same prescription as for $\sG(2n)$. We only modify relation on $\langle h\rangle$ and make it exceptional. Then the defining relation of the odd quadric (\ref{E:relationodd}) can be put into the form  (\ref{E:relationgen}).

\begin{proposition}
Suppose that $\dim(V)$ is odd . Let $A$ be the algebra of algebraic  function on affine non-degenerate quadric in $V^*$. Let $\sG(2n+1)$ be a weak lattice define above. Then  the image of the set of $\sG$-monomials in $A$ is a basis.
\end{proposition}
\begin{proof} The proof repeats the proof in even-dimensional case.
\end{proof}

We are in position to derive a formula for  $A(t)$. The set of  $\sG$-monomials contains  elements of the form
 \[(f^n)^{\alpha_n}\cdots (f^1)^{\alpha_1}g_1^{\beta_1}\cdots g_n^{\beta_n}\] and  \[(f^n)^{\alpha_n}\cdots (f^1)^{\alpha_1}g_1^{\beta_1}\cdots g_n^{\beta_n}h\]    
$\alpha_i,\beta_j\geq0$

These monomials constitute nonintersecting set $E_1$ and $E_2$.

Then \[\begin{split}&A(t)=E_1(t)+E_2(t)\\
&=\frac{1}{(1-t)^{2n}}+ \frac{t}{(1-t)^{2n}}=\frac{1-t^2}{(1-t)^{2n+1}}\end{split}\]

\section{Drinfelds' quasimaps}\label{S:quasimaps}
The combinatorial  approach developed in the last section can be adapted to computation of generating function of Drinfeld's spaces of quasimaps. 

Let us review briefly following  \cite{Braverman} the basics of theory of quasimaps.  
Let $\Sigma$ be an algebraic curve, $X = \mathbf{P}^N$-projective space. In this case points of the space of maps of degree $d$ $\Maps^d(\Sigma, X)$ are classified by  the following data:
\begin{itemize}
\item A line bundle $L$ on $\Sigma$ of degree $-d$.
\item An embedding of vector bundles $L \rightarrow  \O^{N+1}_\Sigma$
\end{itemize}
The reason is that every such embedding defines a one-dimensional subspace in $\mathbb{C}^{N+1}$ for every point $c\in \Sigma$ and thus we get a map $\Sigma\rightarrow \mathbf{P}^N$.

Consider, for example, the case when $\Sigma = \mathbf{P}^1$-projectivization of two-dimensional complex symplectic linear space with a basis $x,y$. In that case $L$ must be isomorphic to the line bundle $\O_{\mathbf{P}^1} (-d)$ (note that such an isomorphism is defined uniquely
up to a scalar) and thus $\Maps^d(\mathbf{P}^1, \mathbf{P}^N )$ becomes an open subset in the projectivization of the vector space $\Hom(\O(-d),\O^{N+1}) \cong \mathbb{C}^{(N+1)}\otimes \bigoplus _{i=0}^d \langle  x^{d-i}y^i\rangle$, i.e. $\Maps^d(\mathbf{P}^1, \mathbf{P}^N ) $
is an open subset of $\mathbf{P}^{(N+1)(d+1)-1}$. The reason that it does not coincide with it is that not every non-zero map $\O (-d)\rightarrow \O^{N+1}$ gives rise a map $\mathbf{P}^1\rightarrow \mathbf{P}^N$ - we
 need to consider only those maps which don't vanish in every fiber.

The above example suggests the following compactification of $\Maps^d(\Sigma,\mathbf{P}^N )$. Namely, we define the space of quasi-maps from $\Sigma$ to $X$ of degree $d$ (denoted by $Q\Maps^d(\Sigma,X)$) to be the scheme classifying the following data:
\begin{enumerate}
\item A line bundle $L$ on $\Sigma$ 
\item  A non-zero map $\kappa : L \rightarrow \O_\Sigma^{N+1}$
\item  Note that $\kappa$ defines an honest map $U\rightarrow \mathbf{P}^N$ where $U$ is an open subset of
$\Sigma$. We require that the image of this map lies in $X$.
\end{enumerate}
For example it is easy to see that if $ X = \mathbf{P}^N$ and $\Sigma =\mathbf{ P}^1$ then $Q\Maps^d(\Sigma,X) \cong \mathbf{P}(\mathbb{C}^{(N+1)}\otimes  \bigoplus _{s=0}^d \langle  x^{d-s}y^s\rangle)$. Let $z^s=x^sy^{d-s}$. Then the affine cone  $CQ\Maps^d(\mathbf{ P}^1, \mathbf{P}^N) $ over projective space $Q\Maps^d(\mathbf{ P}^1,\mathbf{ P}^N)$ is a space of vector-valued polynomials 
\begin{equation}\label{E:polynomial}
\lambda(z)=\sum_{i=0}^d \lambda[s]z^s, \lambda[s]\in \mathbb{C}^{N+1}
\end{equation}
If a submanifold $X\subset \mathbf{P}^N$ is defined by homogeneous equations $r_i(\lambda)=0$ then $CQ\Maps^d(\mathbf{ P}^1, X) $ is defined by equations 
\begin{equation}\label{E:relations}
r_i[s]=\frac{1}{s!}\frac{d^sr_i(\lambda(z))}{dz^s}|_{z=0}=0, s\geq 0
\end{equation}
on coefficients $\lambda[t]$. 
This general theory explains  $\SL(2)$ action on the space  $CQ\Maps^d(\mathbf{ P}^1, X) $ , which is not seen in the description of $CQ\Maps^d(\mathbf{ P}^1, X) $ as a variety  of polynomial maps $\lambda(z)$. 

The spaces of maps of the form \[\lambda(z)=\sum_{i=-r}^{d-r} \lambda[i]z^i\] is isomorphic to the space of maps (\ref{E:polynomial}). The identification, which we will be using  freely,  is made by multiplication of $\lambda(z)$ on $z^r$. 

The main property of the scheme $Q\Maps^d(\Sigma,X)$ is that it possesses a stratification
\[Q\Maps^d(\Sigma,X)=\bigcup_{d'=0}^d\Maps^{d'}(\Sigma,X)\times \Sym^{d-d'}(\Sigma)\]
The corresponding strata in $CQ\Maps^d(\Sigma,X)$ correspond  maps $\lambda(z)$ that can be factored  into   $\tilde{\lambda}(z)f(z)$, where a  scalar polynomial $f(z)$ has degree $d-d'$. The corresponding point in $\Sym^{d-d'}(\Sigma)$ is the zero divisor of $f$. Note that $Q\Maps^d(\Sigma,X)$ is typically singular even for a smooth $X$. Its singularities are located at the described strata.

The Chern class $c_1(\mu)$ of the line bundle $\mu$, that defines an embedding $X\rightarrow \mathbf{ P}^N$ is an element in $H^2(X)$. The degree $d$ of the map $\lambda:\Sigma\rightarrow X$ is the value of the pairing $<c_1(\mu),[\Sigma]>$.  In theory $\Maps^d(\Sigma,X)$ might contain components of  different  dimension.  The (virtual) dimension of $\Maps^d(\Sigma,X)$ can be computed with an aid  of Riemann–Roch formula:
\[\mathrm{vdim}=\dim(X)+<c_1(T_X),[\Sigma]>, <c_1(\mu),[\Sigma]>=d\]  
If it is possible to find two $[\Sigma_1],[\Sigma_2]$ such that $<c_1(T_X),[\Sigma_1]>\neq <c_1(T_X),[\Sigma_2]>$ and $<c_1(\mu),[\Sigma_1]>=<c_1(\mu),[\Sigma_2]>=d$, then reducibility of  $CQ\Maps^d(\Sigma,X)$ is questionable. If  the first Chern class $c_1(T_X)$ of $X$ is very large negative,  then the space $Q\Maps^d(\Sigma,X)$ set-theoretically consists of one stratum $X\times \Sym^{d}(\Sigma)$.  We believe that this unusually degenerate structure of the space  can be a source of anomalies in $\beta\gamma$-system on the cone $CX$.

The spaces of quasimaps to compact homogenous spaces have an interpretation of semi-infinite  closed Schubert cells of a suitable semi-infinite flag variety of an affine Lie algebra $\hat{\g}$ (see e.g. \cite{Braverman}).  Properties of algebras of homogeneous functions on closed Schubert cells in finite dimensional partial flag spaces of finite dimensional semi-simple $\g$ were under intensive scrutiny. Its Koszul property  was established in \cite{Ravi} \cite{Bogvad} \cite{Bezrukavnikov}. Thus Koszul property of its semi-infinite analogs is expected and has been proved for quasimaps to Grassmannian \cite{SottileSturmfels} and to Lagrangian Grassmannian \cite{Ruffo}.


Returning to the case of quasimaps to a quadric we note that  the algebra of functions $AQ$ (defined in the Introduction) carries several gradings.  The most coarse  is defined as the  degree in variables $\lambda^i[k]$. Thus \[\deg \prod_{s}(\lambda^{i_s}[k_s])^{\alpha_s}=\sum_i \alpha_s\]  Then $AQ$ can be decomposed into a direct sum 
\begin{equation}\label{E:grading}
AQ=\bigoplus_{i\geq 0} AQ_i,
\end{equation} where $AQ_i$ is a linear space of elements of degree $i$.
  

Our plan is to compute   generating function (Poincar\'{e} series ) \[AQ_{-N_1}^{N_2}(t)=\sum_{i\geq 0} \dim AQ_{-N_1, i}^{N_2} t^i\] using the example of a quadric as a guide.
We want to stress that our method has been extracted (with suitable modifications) from works \cite{SottileSturmfels} \cite{Ruffo}.

\subsection{Quasi-maps to even-dimensional quadric}


Let us consider  an algebra  of polynomials $\mathbb{C}[f^i[s],g_j[t]]$ with $1\leq i,j\leq n$($n\geq 1$),  $-N_1\leq s,t\leq N_2, N_1,N_2\geq 0.$ 
The algebra \[AQ=AQ_{-N_1}^{N_2}(2n)=\mathbb{C}[f^i[s],g_j[t]]/(r[ -2N_1],\dots,r[ 2N_2])\] is a quotient of a polynomial algebra  by the ideal generated by relations 
\begin{equation}\label{E:relationsloopppp}
r[l]=\sum_{s+t=l}\sum_{i=1}^nf^i[s]g_i[t]=0\quad  -2N_1\leq l\leq 2N_2
\end{equation}
To make a connection with the previous section we note that $f^i[s],g_j[t]$ are Taylor coefficients (after multiplication on $z^{N_1}$) of the coordinate functions of a polynomial map $\lambda(z)$ from $\mathbb{C}^*$ to $V=W+W^*$. 
The relations $r[l]$ are obtained by prescription (\ref{E:relations}) from the equation of the quadric. 

 Dimension of $V$  will be greater or equal to six to the end of the section. 
The generating space $AQ_1$ of $AQ$ (\ref{E:grading})  can be identified with a subspace of $V\otimes \mathbb{C}[z,z^{-1}]$: 
\[\begin{split}
&g_i[t]\rightarrow g_i\otimes z^t\\
&f^i[t]\rightarrow f^i\otimes z^t
\end{split}\]
Operators $R,R_i$ and the elements of the Cartan subalgebra from Section \ref{S:even} act on $AQ_1$.
An additional not everywhere defined on $AQ_1$  operator  $\hat{R}$ will be proved useful.
The space $V\otimes \mathbb{C}[z,z^{-1}]$ is a representation $\rho$ of the Lie algebra $\so_{2n}\otimes \mathbb{C}[z,z^{-1}]$. 
Under identification $\so_{2n}\otimes \mathbb{C}[z,z^{-1}]\cong \Lambda^2(V) \otimes \mathbb{C}[z,z^{-1}]$ operator $\hat{R}$ coincides  with 
\begin{equation}\label{E:rhat}
\rho(\hat{r})=\rho(f^{n-1}\wedge f^{n}\otimes z).
\end{equation}
 It acts by the formula
\[
\hat{R}f^{i}[t]=0 \quad
\hat{R}g_{i}[t]=
\begin{cases}f^n[t+1]  & \text{if $i=n-1$,}\\
-f^{n-1}[t+1]  & \text{if $i=n$,}\\
0 &\text{if $i\neq n-1,n$}\
\end{cases}
\quad 
\]
A reader familiar with theory of affine Lie algebras will immediately  recognize in  $r,r_i,\hat{r}$ positive simple root vectors of  $\hat{\so}_{2n}$.

Mimicking considerations of a quadric we can use operators $R,\hat{R},R_i$ and one-dimensional spaces $\langle g_i[t] \rangle,\langle f^j[t] \rangle$ to define the following affine diagram :


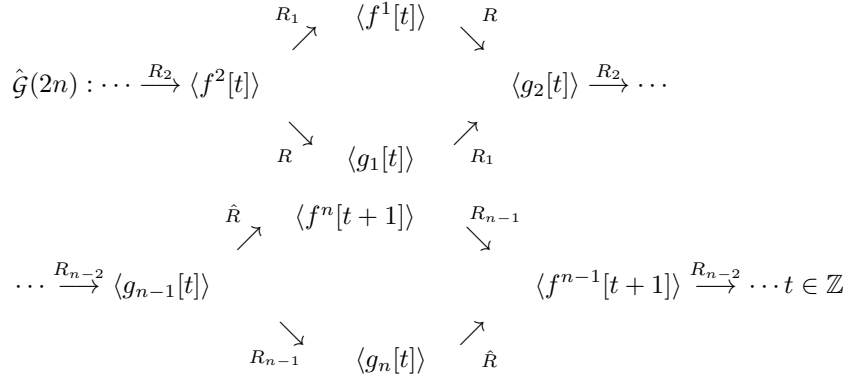
\[\begin{split}
&\hat{\sG}(2n): \cdots\overset{R_2} \longrightarrow  \langle f^{2}[t] \rangle\begin{array}{c}\overset{R_1\quad \text{ }}\nearrow  \begin{array}{c} \langle f^1[t] \rangle\\ \text{ } \end{array} \overset{\text{ }\quad R}\searrow \\ \underset{R \quad  \text{ }}\searrow   \begin{array}{c} \text{ }\\\langle  g_1[t] \rangle \end{array} \underset{\text{ }\quad R_1}\nearrow \end{array}\langle g_2[t] \rangle \overset{R_2}\longrightarrow \cdots \\
&\cdots\overset{R_{n-2}}  \longrightarrow  \langle g_{n-1}[t] \rangle\begin{array}{c}\overset{\hat{R}\quad \text{ }}\nearrow  \begin{array}{c} \langle f^{n}[t+1] \rangle \\ \text{ } \end{array} \overset{\text{ }\quad R_{n-1}}\searrow \\ \underset{R_{n-1} \quad  \text{ }}\searrow   \begin{array}{c} \text{ }\\\langle g_n[t] \rangle \end{array} \underset{\text{ }\quad \hat{R}}\nearrow \end{array} \langle f^{n-1}[t+1] \rangle \overset{R_{n-2}}\longrightarrow\cdots t\in \mathbb{Z}
\end{split}
 \]

This diagram defines an infinite poset  \[\hat{\sG}=\hat{\sG}(2n)=\{e_{\alpha}[t]|e_{\alpha}\in \sG(2n),t\in \mathbb{Z}\}.\]
 A finite sub-poset  $\sG_{-N_1}^{N_2}$ labels  subspaces of the first component of  $AQ_{-N_1}^{N_2}$. Alternatively, poset $\sG_{-N_1}^{N_2}$ is an interval
\begin{equation}\label{E:interval}
 G_{-N_1}^{N_2}=[\langle f^n[-N_1]\rangle ,\langle g_n[N_2]\rangle ]\overset{\ddef}{=}\{\gamma|\langle f^n[-N_1]\rangle\leq \gamma \leq \langle g_n[N_2]\rangle\}
 \end{equation}

The equations (\ref{E:relationsloopppp}) can be rewritten in a form similar to (\ref{E:relationzero})
\begin{equation}\label{E:relationN}
\begin{split}
&g_1[t]f^1[t]=-g_2[t]f^2[t]-\sum_{ 3\leq i\leq n} g_i[t]f^i[t]-\sum_{s\neq 0}\sum_{i=1}^n g_i[t+s]f^i[t-s]\\
&g_n[t]f^n[t+1]=-g_{n-1}[t]f^{n-1}[t+1]-\sum_{1\leq i\leq n-2} g_i[t]f^i[t+1]-\\
&-\sum_{s\neq 0}\sum_{i=1}^n g_i[t+s]f^i[t+1-s]
\end{split}
\end{equation}
Index  $s$ in the summations belongs to a maximal subset of integers such that all monomials are elements of $AQ_{-N_1}^{N_2}$ 

If we label uniformly generators $e_{\alpha}$ of $AQ_{-N_1}^{N_2}$ by elements of  $\sG_{-N_1}^{N_2}$, then relations (\ref{E:relationN}) would have a form (\ref{E:relationgen}).

\begin{proposition}
Suppose $n\geq 3$. 
Introduce a total order on $\hat{\sG}(2n)$, which is a refinement of the partial order.
We set $\langle f^1[t] \rangle$ to be grater then $\langle g_1[t] \rangle$ and $\langle f^n[t+1]\rangle$ to be grater then $\langle g_n[t]\rangle$. Then the generator (\ref{E:relationzero}) of the ideal of relations defines a Gr\"{o}bner basis of the ideal with respect to the degree-lexicographic order on monomials.
\end{proposition}
\begin{proof}
A short introduction Gr\"{o}bner bases technique is given in Appendix \ref{S:Grobner}. Here we follow notations of that section.

We need to compute appropriate $S$-polynomials. Note that the leading monomials of relations $\{r[l]\}$ (\ref{E:relationN}) are  $g_1[t]f^1[t]$ and $g_n[t]f^n[t+1]$. If $n\geq 3$,  then these monomials are relatively prime in the semigroup generated by $g_i[s],f^j[t]$. By Lemma (\ref{L:prime}) the set of $S$-polynomials, that we have to compute in order to find $Gr_{1}(I(r[l]))$, is empty. Thus $Gr(I(r[l]))=\{r[l]\}$.

\end{proof}

One can prove the following theorem
\begin{theorem}\label{T:basisloop}
The image of the set of $\sG_{-N_1}^{N_2}(2n)$-monomials in $QA_{-N_1}^{N_2}(2n)$ is a basis. Thus $QA_{-N_1}^{N_2}(2n)$ is an algebra with a straightened law. 
\end{theorem}
\begin{proof}
We can interpret results of the previous proposition as follows. The semigroup ideal $T(\{r[l]\})$ coincides with the set of elements  divisible by $g_1[t]f^1[t]$, $g_n[t]f^n[t+1]$. These are precisely  non-$\sG_{-N_1}^{N_2}$-monomials.  The complement of non-$\sG_{-N_1}^{N_2}(2n)$-monomials in the semigroup of monomials is the set of   $\sG_{-N_1}^{N_2}(2n)$-monomials. By Theorem  \ref{T:basic1} item  \ref{I:one} the image of this set in $QA_{-N_1}^{N_2}(2n)$  is a basis.
\end{proof}


\subsubsection{Digression about quadratic and Koszul algebras}\label{S:Koszul}

Recall, that a graded (not necessarily commutative) algebra \ $A=\bigoplus_{n\geq 0} A_{n}$ \ is
 a quadratic if \ $A_{0}=\mathbb{C}, $ \ \ $W=A_{1}$ \
generates \ $A$ \ and all relations follow from quadratic relations \ $\sum\limits_{i, j}r_{ij}^{k}x^{i}x^{j}=0$ \ where \ $x^{1}, \ldots , x^{\dim W}
$ \ is a basis of \ $W=A_{1}.$ \ The space of quadratic relations $R.$  is the
subspace of \ $W\otimes W$ \ spanned by \ $r_{ij}^{1}, r_{ij}^{2}, \ldots $)
.  Then \ $A_{2}=W\otimes W/R$ \ and \ 
$A$ \ is a quotient of free algebra (tensor algebra) \ $\bigoplus_{n\geq 0} W^{\otimes n}$ \ with
respect to the ideal generated by \ $R.$ \ The dual quadratic algebra \ $%
A^{!}$ \ is defined as a quotient  \ $\bigoplus_{n\geq 0} \left(W^{\ast }\right) ^{\otimes n}/I(R^{\perp })$. The ideal $I(R^{\perp })$ of the free algebra $\bigoplus_{n\geq 0} \left(W^{\ast }\right) ^{\otimes n}$ is 
\ generated  by \ $R^{\perp }\subset W^{\ast }\otimes W^{\ast }$ \ (here \ $%
R^{\perp }$ \ stands for the subspace of \ $W^{\ast }\otimes W^{\ast
}=\left(W\otimes W\right) ^{\ast }$ \ that is orthogonal to \ $R\subset
W\otimes W$). 

The case of commutative $A$ has its specifics. Quadratic relations of algebra $A$ contains commutators $x^ix^j-x^j x^i$. An easy computation shows that $R^{\perp}$ contains only anti-commutators. We can interpret $A^{!}$ as a universal enveloping of some Lie algebra $L$. In the following we will refer to $L$ as to Koszul dual to commutative $A$. 

Let $\g$ be a super Lie algebra with structure constants $f_{\alpha\beta}^{\gamma}$ in basis $b_{\alpha}$. We equip  algebra  of functions on generators $c^{\alpha}$ dual to $b_{\alpha}$, having the reverse parity an operator $d$. It acts  by the formula \[d=f_{\alpha\beta}^{\gamma}c^{\alpha}c^{\beta}\frac{\sd}{\sd c^{\gamma}}=f_{\alpha\beta}^{\gamma}c^{\alpha}c^{\beta}b_{\gamma}\] and satisfies $d^2=0$.

This formula reminiscences BRST differential in the gauge theory. In mathematical literature (see e.g. \cite{Chevalley} ) the complex $(\mathbb{C}[c^{\alpha}],d)$  is known as Cartan-Chevalley  complex of Lie super-algebra $\g$ and is usually denoted by $C(\g)=C^{\bullet}(\g)$.

We will give  a definition of Koszul algebra that is most suited to our purposes. Other definitions and their equivalence are discussed in \cite{PP}.
\begin{definition}
A  quadratic commutative algebra $A$ is called Koszul if the cohomology $H(L)$ of the complex $C(L)$ for Koszul dual $L$ is equal to $A$.
\end{definition}

\begin{theorem}
The algebra $QA_{-N_1}^{N_2}(2n)$ ($n\geq 3$) is Koszul .
\end{theorem}
\begin{proof}
Follows from straightened law property of the algebra $QA_{-N_1}^{N_2}(2n)$ (see \cite{Kempf}).
\end{proof}

\subsection{Quasi-maps to odd-dimensional quadric}
The relevant  algebra  of polynomials  is $\mathbb{C}[f^i[s],g_j[t],h[k]]$ with $1\leq i,j\leq n$($n\geq 1$),  $-N_1\leq s,t,k\leq N_2, N_1,N_2\geq 0.$ 
The algebra 
\begin{equation}\label{E:algebraodd}
AQ=AQ_{-N_1}^{N_2}(2n+1)=\mathbb{C}[f^i[s],g_j[t],h[k]]/(r[ -2N_1],\dots,r[ 2N_2])
\end{equation} is a quotient by the   relation ideal
\begin{equation}\label{E:relationsodd}r[l]=\sum_{s+t=l}\left(\sum_{i=1}^nf^i[s]g_i[t]\right)+h[s]h[t]=0\quad  -2N_1\leq l\leq 2N_2\end{equation}
The element $\hat{r}$ (see equation \ref{E:rhat}) gets modified as follows
\[\hat{r}=f^{n-1}\wedge f^{n}\otimes z\]
This enables us to find the corresponding periodic (affine) Hasse diagram ($n\geq 2$)
\[\begin{split}
& \hat{\sG}(2n+1): \cdots\overset{R_1} \longrightarrow \langle f^{1}[t] \rangle \overset{R}\longrightarrow  \langle h[t]\rangle  \overset{R}\longrightarrow \langle g_1[t] \rangle \overset{R_1}\longrightarrow \cdots \\
&\cdots\overset{R_{n-2}}  \longrightarrow \langle g_{n-1}[t] \rangle\begin{array}{c}\overset{\hat{R}\quad \text{ }}\nearrow  \begin{array}{c}\langle f^{n}[t+1] \rangle\\ \text{ } \end{array} \overset{\text{ }\quad R_{n-1}}\searrow \\ \underset{R_{n-1} \quad  \text{ }}\searrow   \begin{array}{c} \text{ }\\\langle g_n[t] \rangle \end{array} \underset{\text{ }\quad \hat{R}}\nearrow \end{array}\langle f^{n-1}[t+1] \rangle \overset{R_{n-2}}\longrightarrow\cdots t\in \mathbb{Z}
\end{split}
 \]

 We can construct a weak lattice (see Definition \ref{D:wlattice}) from the above diagram \[\hat{\sG}(2n+1)=\{e_{\alpha}[t]|e_{\alpha}\in \sG(2n+1),t\in \mathbb{Z}\}\] Construction is the same as in the even case, except that non-reflexive elements  $h[t]$  are present. 
 We define $\hat{\sG}_{-N_1}^{N_2}(2n+1)$ by the formula (\ref{E:interval}).
 
 Relations  (\ref{E:relationsodd}) can be put in the form 
 \begin{equation}\label{E:relationoddN}
\begin{split}
&(h[t])^2=-g_1[t]f^1[t]-\\
&-\sum_{i=2}^n g_i[t]f^i[t]\\
&-\sum_{s\neq 0}\left( h[t+s]h[t-s] +\sum_{i=1}^n g_i[t+s]f^i[t-s]\right)\\
&g_n[t]f^n[t+1]=-g_{n-1}[t]f^{n-1}[t+1]- \\
&-\sum_{i=1}^{n-2} g_i[t]f^i[t+1]-\\
&-\sum_{s\neq 0} \left( h[t+1+s]h[t-s] +\sum_{i=1}^n g_i[t+s]f^i[t+1-s]\right)
\end{split}
\end{equation}

\begin{proposition}
Suppose $n\geq 2$. 
Introduce a total order on $\hat{\sG}_{-N_1}^{N_2}(2n+1)$, which is a refinement of the partial order.
We set $\langle f^1[t] \rangle <\langle h[t]\rangle < \langle g_1[t]\rangle$ and $\langle f^n[t+1]\rangle$ to be grater then $\langle g_n[t]\rangle$. Then the generators (\ref{E:relationoddN}) of the ideal of relations defines a Gr\"{o}bner basis of the ideal with respect to the degree-lexicographic order on monomials.
\end{proposition}
\begin{proof} It is  similar to the proof of Proposition \ref{T:basisloop}.

We need to compute appropriate $S$-polynomials. Note that the leading monomials of relations $\{r[l]\}$ (\ref{E:relationN}) are  $h[t]^2$ and $g_n[t]f^n[t+1]$. If $n\geq 1$,  then these monomials are relatively prime in the semigroup generated by $g_i[s],f^j[t],h[k]$. By Lemma (\ref{L:prime}) the set of $S$-polynomials, that we have to compute in order to find $Gr_{1}(I(r[l]))$, is empty. Thus $Gr(I(r[l]))=\{r[l]\}$.

\end{proof}

As a corollary we get the following 
\begin{theorem}\label{T:basisloopodd}
The image of the set of $\hat{\sG}_{-N_1}^{N_2}(2n+1)$-monomials in $QA_{-N_1}^{N_2}(2n+1)$ is a basis.
\end{theorem}

\begin{theorem}
 $QA_{-N_1}^{N_2}(2n+1)$ is a Koszul algebra.
\end{theorem}
\begin{proof}
Formally  $QA_{-N_1}^{N_2}(2n+1)$ is not an algebra with straightened law. However the proof of \cite{Kempf} can be easily adapted to our case. 

Alternatively the basis of  $\hat{\sG}_{-N_1}^{N_2}(2n+1)$ monomials is  Poincar\'{e}-Birkhoff-Witt basis (see \cite{PP} for definition).
Priddy in \cite{Priddy} have proven that any  quadratic PBW-algebra is Koszul. 
\end{proof}
\section{Generating functions}\label{S:generating}
Theorem \ref{T:basisloop} enables us to find various generating function related to algebra   $QA_{-N_1}^{N_2}$, besides  $QA_{-N_1}^{N_2}(t)$.
The group of obvious $\SO(n)$-symmetries of algebra  $QA_{-N_1}^{N_2}$ can be extended to $\SO(n)\times \SL(2)$. The action of $\SL(2)$ factor is explained in Section \ref{S:quasimaps}. 


It makes sense to define a generating function of  $\SO(n)\times \SL(2)$-characters $QA(g_1,g_2,t)$. This formal function is completely determined by its restriction $QA(z_1,\dots,z_{ [\frac{n}{2}]},q,t)$ on the maximal complex torus $\mathbb{C}^{\times [\frac{n}{2}]}\times \mathbb{C}^{\times} \subset \SO(n)\times \SL(2)$. 

For simplicity of exposition we let $z_i=1$ . Theorem \ref{T:basisloop} reduces computation of $QA(q,t)$  to  combinatorics. We can associate to any $k$-tuple  \[m=(e_{\alpha_1}[l_1],\dots, e_{\alpha_k}[l_k])\in \bigcup_{k\geq 1} \hat{\sG}^{\times k}/\Sigma_k\]  a local weight \[w(m)=\prod_{s=1}^k q^{l_s}\] and  a function $\deg(m)=k$. Let $E$ be some set of $m$. We define a generating function \[\hat{E}(t)=\sum_{m\in E}w(m)t^{\deg(m)}.\]
The direct corollary of Theorem \ref{T:basisloop} is that \[QA_{-N_1}^{N_2}(q,t)=T_{-N_1}^{N_2}(q,t),\] where $T_{-N_1}^{N_2}$ is a subset of  $\bigcup_{k\geq 1} \hat{\sG}_{-N_1}^{N_2\ \times k}/\Sigma_k$ which consists of all $\hat{\sG}_{-N_1}^{N_2}$-elements.


{\bf Even-dimensional quadric}

We will be interested in the intervals \[\Delta_l=[f^{n-2}[l], f^{n-1}[l+1]].\] Then $\bigcup_{l}\Delta_l=\hat{\sG}.$
Smaller subintervals will also be used:
\[\Delta'_l=[f^{n}[l], f^{n-1}[l]]\]
\[\Delta''_l=[f^{n-2}[l], g_{n}[l]]\]
The set $\hat{\sG}_{-N_1}^{N_2}$ is a union \[\hat{\sG}_{-N_1}^{N_2}=\Delta'_{-N_1} \cup \left( \bigcup_{l=-N_1}^{N_2-1} \Delta_t \right)\cup \Delta''_{N_2}\]

 The last equality implies that \[T_{-N_1}^{N_2}=T(\Delta'_{-N_1}) \times  \left( \prod_{l=-N_1}^{N_2-1} T(\Delta_l) \right)\times  T(\Delta''_{N_2})\]
 This allows us to write a formula for generating function 
 \[T_{-N_1}^{N_2}(q,t)=T(\Delta'_{-N_1})(q,t) T(\Delta''_{N_2})(q,t) \prod_{l=-N_1}^{N_2-1} T(\Delta_l)(q,t)  \]
 It is not hard to find formulas for local factors:
 
\[T(\Delta'_{l})(q,t)=\frac{1}{(1-q^lt)^2}\]
The formula for $T(\Delta''_{l})(q,t)$ and $T([f^{n-2}[l], g_{2}[l]])(q,t)$ can be obtained along the same lines as for  $A(t)$:

\[T(\Delta''_{l})(q,t)=\frac{1-q^{2l}t^2}{(1-q^lt)^{2n-2}}\]
\[T([f^{n-2}[l], g_{2}[l]])(q,t)=\frac{1-q^{2l}t^2}{(1-q^lt)^{2n-(2+n-2)}}=\frac{1-q^{2l}t^2}{(1-q^lt)^{n}}\]

Finally $\Delta_l=[f^{n-2}[l], g_{2}[l]]\cup[g_{3}[l], f^{n-1}[l+1]]$. Using inclusion-exclusion principle we get 

 \[\begin{split}
 &T([g_{3}[l], f^{n-1}[l+1]])(q,t)=T([g_{3}[l], g_{n-2}[l+1]])(q,t) T([g_{n-1}[l], f^{n-1}[l+1]])(q,t)=\\
 &\frac{1}{(1-q^lt)^{n-4}}\left(\frac{1}{(1-q^lt)(1-q^lt)(1-q^{l+1}t)}+\frac{1}{(1-q^lt)(1-q^{l+1}t)(1-q^{l+1}t)}-\frac{1}{(1-q^lt)(1-q^{l+1}t)}\right)=\\
 &=\frac{1}{(1-q^lt)^{n-4}}\frac{1-q^{2l+1}t^2}{(1-q^lt)^2 (1-q^{l+1}t)^2}
 \end{split}\]

Finally 
 \[\begin{split}&T_{-N_1}^{N_2}(q,t)=  \frac{1-q^{2N_2}t^2}{(1-q^{-N_1}t)^2(1-q^{N_2}t)^{2n-2}}  \prod_{l=-N_1}^{N_2-1} \frac{(1-q^{2l}t^2)(1-q^{2l+1}t^2)}{(1-q^lt)^{n}(1-q^lt)^{n-4}(1-q^lt)^2 (1-q^{l+1}t)^2}=\\
&=\frac{\prod_{l=-2N_1}^{2N_2}(1-q^lt^2)}{\prod_{l=-N_1}^{N_2}(1-q^lt)^{2n}} \end{split}  \]

{\bf Odd-dimensional quadric}

The proof of the formula 

\[T_{-N_1}^{N_2}(q,t)=\frac{\prod_{l=-2N_1}^{2N_2}(1-q^lt^2)}{\prod_{l=-N_1}^{N_2}(1-q^lt)^{2n+1}},n\geq 1\]

is similar to even-dimensional case and is left to the interested reader as an exercise.
We have proved a 
\begin{proposition}
The generating function $AQ_{-N_1}^{N_2}(n)(q,t)$ , $n\geq 5,N_1,N_2\geq 0$ is 
\[\frac{\prod_{l=-2N_1}^{2N_2}(1-q^lt^2)}{\prod_{l=-N_1}^{N_2}(1-q^lt)^{n}}\]
\end{proposition}
We will treat the remaining $n=2,3,4$ cases in the appendix
\section{The mini-BRST resolution}

 A special interest is the value of $T_{-N_1}^{N_2}(q,t)$ at $q=1$:
 
 \begin{equation}\label{E:genfunct}
 T_{-N_1}^{N_2}(1,t)=\frac{(1-t^2)^{2(N_1+N_2)+1}}{(1-t)^{n(N_1+N_2+1)}}
 \end{equation}
 
We know that algebras  $AQ=AQ_{-N_1}^{N_2}(n)$ $n\geq 5,N_1,N_2\geq 0$  Koszul . The  Koszul dual algebra $\left(AQ_{N_1}^{N_2 }(n)\right)^!$ (the universal enveloping algebra of a graded Lie algebra $L$, described in Section \ref{S:Koszul}) has (see \cite{PP}) the generating function 
\begin{equation}\label{E:dual}
1/T_{N_1}^{N_2}(1,-t)=\frac{(1+t)^{n(N_1+N_2+1)}}{(1-t^2)^{2(N_1+N_2)+1}}.
\end{equation}
For this  Lie algebra we introduce a special notation:
\[LH=LH_{-N_1}^{N_2}(n)=\bigoplus_{k\geq 0} LH_k\] It has  finite-dimensional graded components $LH_k$ (we suppress  dependence on $n,N_1,N_2$).

The formal power series  satisfy  \[AQ(-t)^{-1}=\frac{\prod_{k\geq0} (1+t^{2k+1})^{\dim LH_{2k+1}}} {\prod_{k\geq1} (1-t^{2k})^{\dim LH_{2k}}}\] The last identity represents 
numeric version of Poincar\'{e}-Birkhoff-Witt  theorem. What is important is that this identity allows unambiguously recover dimensions (and $\SO(n)\times\SL(2)$ characters) of $LH_k$. 
Analyzing  formula (\ref{E:dual}) we conclude that 
$L_1$ has dimension 
$n(N_1+N_2+1)$ ,  $L_2$ has dimension  $2(N_1+N_2)+1$ and all other graded components vanish.

 We will   describe Lie algebra $L'H_{-N_1}^{N_2}$. Let $H$ be an odd  graded Heisenberg algebra. It is generated by linear space $V$ in degree one. It also contains a central element $b$ in degree two, which  satisfies 
 \begin{equation}\label{E:commutator}
 [v_1,v_2]=(v_1,v_2)b\in H_1.
 \end{equation} The bracket stands for the graded commutator (in our case anti-commutator).  We define the loop algebra $\hat{H}=H\otimes \mathbb{C}[z,z^{-1}]$.  The algebra $\hat{H}$ contains a subalgebra \[L'H_{-N_1}^{N_2}=\bigoplus_{l=-N_1}^{N_2} V\otimes z^l+ \bigoplus_{l=-2N_1}^{2N_2} \langle b \rangle \otimes z^l\]
\begin{proposition}
The universal enveloping $U(L'H_{-N_1}^{N_2}(n))$ is isomorphic to $\left(AQ_{-N_1}^{N_2 }(n)\right)^!$
\end{proposition}

\begin{proof}

The space of generators $LH=LH_{-N_1}^{N_2}(n)$ is dual to the corresponding space of $AQ_{-N_1}^{N_2 }(n)$ which is equal to $AQ_1=\bigoplus_{l=-N_1}^{N_2} V\otimes z^l$. Thus $LH_1=L'H_1=\bigoplus_{l=-N_1}^{N_2} V\otimes z^l$ (we identify $V$ with its dual by means of inner product). It is known (see \cite{PP}) that  the space of elements of degree two $LH_2$ is dual to the space of relations of $AQ_{-N_1}^{N_2 }(n)$ (\ref{E:reltrivial}). The later can be identified with $\bigoplus_{l=-2N_1}^{2N_2} \langle b \rangle \otimes z^l$.

The Lie algebra $LH$ is generated by $LH_1$ (this is equivalent to Koszul property of $AQ_{-N_1}^{N_2 }(n)$ \cite{PP}). The only possible nonzero  $\SO(n)\times \SL(2)$-invariant graded commutator \[\left(\bigoplus_{l=-N_1}^{N_2} V\otimes z^l\right)^{\otimes 2}\rightarrow  \bigoplus_{l=-2N_1}^{2N_2} \langle b \rangle \otimes z^l\] is given (up to a constant factor ) by  the formula (\ref{E:commutator}).
\end{proof}


\begin{pf}{\ref{T:main}}

Let us now put all our ingredients together. The algebra $C=\mathbb{C}[\lambda^i[s]]\otimes \Lambda[c[k]]$ $i=1,\dots,n, s=-N_1,\dots,N_2, k=-2N_1,\dots,2N_2$ has two interpretations. In the first it is a truncated version of BRST complex from Section  \ref{S:introduction}. In the second interpretation it is a Cartan-Chevalley complex of Lie algebra $LH_{-N_1}^{N_2}(n)$ $n\geq 5$, whose cohomology is equal to $AQ_{-N_1}^{N_2}$. 
The quasi-isomorphism $ \zeta:C\rightarrow AQ$ is defined on generators by the formula \[\zeta(\lambda^i[s])=\lambda^i[s], \zeta(c[s])=0.\]
The exceptional cases $n=3,4$ are handled in Appendix.
\end{pf}

As a side remark we note  that algebraic variety $Spec(AQ_{-N_1}^{N_2}(n))$ is a complete intersection for the above range of parameters $N_1,N_2,n$.


\section{The Hilbert space}\label{S:Hilbert}
One  possible  approach to the  problem of defining  Hilbert space in $\beta\gamma$-system on a quadric \cite{AA} utilizes a technique of semi-infinite cohomology. It coincides with quantum mini-BRST complex.
Let $W$ be a vector space, decomposed into a  direct sum $W=W_+\oplus W_-$.
Any operator $a:W\rightarrow W$ has a block-form $a=\left( \begin{smallmatrix}   a_{++} & a_{+-} \\ a_{-+} & a_{--} \end{smallmatrix}\right)$.

Let  $\g = \bigoplus_{i\in \mathbb{Z}}\g_i$ be a graded superalgebra. Define linear spaces \[U_-= \bigoplus_{i\leq 0}\g_i, U_+= \bigoplus_{i\geq 1}\g_i\]  A linear space $W_-$ is commensurable  with $U_-$ if  projection on $U_-$ has finite dimensional kernel and co-kernel.  There are also $U_+$ commensurable liner spaces. 

The following is a slight modification of a definition taken from \cite{Voronov}.
\begin{definition}  Let $\g$ be a (super) Lie algebra over the field of complex numbers. We say that $\g$ is provided with a semi-infinite structure if the following
data 1 and 2 are given and condition 3 is satisfied:
\begin{enumerate}
\item  a $\mathbb{Z}$-grading $\g = \bigoplus_{i\in \mathbb{Z}}\g_i$ on $\g$ (that implies $[\g_i,\g_j ] \subset \g_{i+j}$ for each $i,j$), such that
$dim \g_i < \infty$ for each $i$.
\item  a $1$-cochain $\beta$ on $\g$, such that its coboundary $\sd \beta = \mathrm{ad}^* \phi$, $\phi$ being the $2$-cocycle
 on $\g$ .
 The two-cocycle  $\phi(a,b)$ is defined b the formula\[\phi(a,b)=\tr(ad(b)_{-+}ad(a)_{+-}-ad(a)_{-+}ad(b)_{+-})\]
 Here $tr$ stands for the (super) trace. Operators $a_{\pm\pm},b_{\pm\pm}$ are components of adjoint operators computed with respect to decomposition \[\g=
 W_-+W_+,\] where $W_-( W_+)$ are $U_-(U_+)$ commensurable liner spaces
\item  the cochain $\beta$ vanishes on $\g_i$ for all $i\neq0$
\end{enumerate}
\end{definition}

The Lie algebra $LH=LH_{-N_1}^{N_2}$ is graded by the powers in $z$. Also it has a $U_-,U_+$ commensurable polarization 
\[W_+=LH_{\geq 1}=\bigoplus_{l=1}^{N_2} V\otimes z^l+ \bigoplus_{l=2}^{2N_2} \langle b \rangle \otimes z^l\]
\[W_-=LH_{\leq 0}=\bigoplus_{l=-N_1}^{0} V\otimes z^l+ \bigoplus_{l=-2N_1}^{1} \langle b \rangle \otimes z^l\]
The space $LH+LH^*$ has a canonical symmetric inner product, which enables us to define a Clifford algebra $Cl(LH+LH^*)$. 
The space $LH+LH^*$ contains an isotropic subspace $LH_{\geq 1}+(LH_{\geq 1})^{\perp}\cong LH_{\geq 1}+(LH_{\leq 0})^*$ (we use the graded dual). 
We use it 
to construct the Fock space.  The Fock space is the Clifford module $\Lambda^{\frac{\infty}{2}}$ generated by a vacuum vector $\omega$ that is annihilated by $LH_{\geq 1}+(LH_{\leq 0})^*$. It   is called  a  module  of semi-infinite forms. Let $\psi_{\alpha}$ be a basis in $LH$ and $\psi^{* \alpha}$ be the dual basis $LH^*$. Denote by $c_{\alpha\beta}^{\gamma}$ the structure constants of $LH$ :\[[\psi_{\alpha},\psi_{\alpha}]=c_{\alpha\beta}^{\gamma}\psi_{\gamma}\] 
\begin{proposition}
The element $d=c_{\alpha\beta}^{\gamma}\psi^{* \alpha}\psi^{* \beta}\psi_{\gamma}$  in Clifford algebra satisfies $d^2=0$. 
\end{proposition}
\begin{proof}
By   general theory (see \cite{Voronov}) the operator $d^2$, acting in $\Lambda^{\frac{\infty}{2}}$,  coincides with the action of operator of multiplication on  $\phi(a,b)$, thought of as an element of Clifford algebra. Note that $LH$ is a nilpotent Lie  algebra with all double commutators $[a[a',a'']]$ equal to zero. Furthermore $[a,ad(a')_{\pm\pm}(a'')]=0$. From this we conclude that $ad(a')_{-+}ad(a)_{+-}-ad(a)_{-+}ad(a')_{+-}=0$ and the cocycle vanishes identically.
\end{proof}

The operator $d$ defines a differential in $\Lambda^{\frac{\infty}{2}}$.
Being a cyclic module over $Cl(LH+LH^*)$  the space $\Lambda^{\frac{\infty}{2}}$  has a $\mathbb{Z}$ grading \[\Lambda^{\frac{\infty}{2}}=\bigoplus_{i\in \mathbb{Z}} \Lambda^{\frac{\infty}{2}+i}.\] Creation operators (operators of multiplication on $a\in LH^* $ ) increase this grading by one, annihilation operators (operators of multiplication on $a'\in LH $) decrease it by one. The space $\Lambda^{\frac{\infty}{2}}$ is isomorphic to the tensor product \[\Lambda(LH_{\geq 1})^*\otimes \Lambda(LH_{\leq 0})\otimes \omega\] of graded exterior powers. 

The generating function of euler characters 
 of this complex is \[\frac{\prod_{l=0}^{2N_2}(1-q^lt^2)}{\prod_{l=0}^{N_2}(1-q^lt)^{n}}           \frac{\prod_{l=1}^{2N_1}(1-q^lt^{-2})}{\prod_{l=1}^{N_1}(1-q^lt^{-1})^{n}}   \]
with an assumption that the grading of vacuum is zero.

The limiting function $N_1,N_2\rightarrow \infty$ \[Z(q,t)=  \frac{1-t^2}{(1-t)^{n}}   \prod_{l=1}^{\infty} \frac{(1-q^lt^2)(1-q^lt^{-2})}{(1-q^lt^{})^{n}(1-q^lt^{-1})^{n}}   \] satisfies  equations
\[Z(q,t^{-1})=-(-t)^{n-2}Z(q,t)\] \[Z(q,qt)=(-1)^{n}t^{n-4}q^{-1}Z(q,t)\]
which can be established by simple manipulation with the product.
Also 
\[Z(q,qt^{-1})=-t^{2}q^{-1}Z(q,t),\]
which is a corollary of the previous equations. The same formulas has been found in \cite{AA}.

One of the results of Voronov \cite{Voronov} is that there is a spectral sequence, that converges to semi-infinite cohomology.

\begin{theorem}\label{T:voronov} For a Lie algebra $\g = \b \oplus \n$ with a semi-infinite structure and a
$\g$-module $M \in \O$, there exists a spectral sequence
\begin{equation}\label{E:spectral}
\{E_r^{p,q}, d^{p,q}_r:E_r^{p,q}  \rightarrow E_r^{p+r,q-r+1}| p\leq 1, q\geq 0,r\geq 0\}
\end{equation}
with the following properties:
\begin{enumerate}
\item $E_ 1^{p,q}= H^q(\n, \Lambda^{\frac{\infty}{2}+p}(\g/\n)^*\otimes M)$.
\item  If $\n$ is an ideal, then $E_2^{p,q}= H^{\frac{\infty}{2}+p}(\g/\n, H^q(\n, M))= H_{-p}(\g/\n, H^q(\n, M) \otimes 
 \Lambda^{\frac{\infty}{2}+0}(\g/\n)^*)$. As usual $H_i(\g)$ stands for Lie algebra homology .
\item  $E_{\infty}^{p,q} = \mathrm{gr}^pH^{\frac{\infty}{2}+p+q}(\g, M)$.
\item  the differentials $d_r$, induce a sequence of epimorphisms $E^{p,q}_r\rightarrow E^{p,q}_r\rightarrow \cdots  $ for $r$ large enough, so that
$\underset{\rightarrow }\lim E^{p,q}_r = E^{p,q}_{\infty}$ .
\end{enumerate}
\end{theorem}
The theorem admits an obvious modification for super-algebras.
The one dimensional  $LH_{\geq 1}$ module $\Lambda^{\frac{\infty}{2}+0}(\g/\n)^*=\Lambda^{\frac{\infty}{2}+0}(LH/LH_{\geq 1})^*$ (with respect to adjoint action) is trivial.  Hence \[\Lambda^{\frac{\infty}{2}-p}(LH/LH_{\geq 1})^*\cong \Lambda^{p}(LH/LH_{\geq 1})\]
  The sequence adapted to our needs has a form:
\[H^i(LH_{\geq 1},\Lambda^j(LH/LH_{\geq 1}))\Rightarrow H^{\frac{\infty}{2}+i-j}(LH,\mathbb{C})\]

There is no convergence issues with this sequence because of the grading induced by $\mathbb{C}^*$ action $z\rightarrow az$ on the loop parameter. The sequence breaks down into a direct sum 
of finite-dimensional complexes according to this grading, for which convergence issues are vacuous.  

By construction the space $\Lambda^{\frac{\infty}{2}}(LH)$ is a module over the Koszul complex  \[C=C(LH)=\Lambda(LH)^*.\] Here we interpret it as a complex  of cohomological chains. We already know (Theorem \ref{T:main}) that the    homomorphism \[C=C(LH)\overset{\zeta}\rightarrow AQ^{N_2}_{-N_1}\] is a quasi-isomorphism. There is a restriction  homomorphism $AQ^{N_2}_{-N_1}\rightarrow AQ_1^{N_2}.$We define a complex \[\Lambda^{\frac{\infty}{2}}(LH)_{red}\overset{\ddef}=\Lambda^{\frac{\infty}{2}}(LH)\underset{C(LH)}{\otimes} AQ_1^{N_2}\]
and  a  homomorphism of $C(LH)$-modules  \[red:\Lambda^{\frac{\infty}{2}}(LH)\rightarrow \Lambda^{\frac{\infty}{2}}(LH)_{red},\] defined by the formula $a\rightarrow a\otimes 1$

\begin{proposition}
The map $red$ is a quasi-isomorphism.
\end{proposition}
\begin{proof}
Filtration on $\Lambda^{\frac{\infty}{2}}(LH)$ that defines spectral sequence  \ref{E:spectral} (see \cite{Voronov} for details) induces a filtration on $\Lambda^{\frac{\infty}{2}}(LH)_{red}$  The map $red$ defines  a map of spectral sequences. The map of $E_1$ terms  
\[red:H(LH_{\geq 1},\Lambda^j(LH/LH_{\geq 1}))\rightarrow H(\Lambda^j(LH/LH_{\geq 1})\otimes AQ_1^{N_2})\]
It is a quasi-isomorphism because $U(LH_{\geq 1})\cong U(LH^{N_2-1}_{0})$ is a Koszul algebra. The  cohomology of any graded $LH_{\geq 1}$ module $M$ (which in our case is  $\Lambda^j(LH/LH_{\geq 1})$) can be computed with the complex $M\otimes U(LH_{\geq 1})^!\cong M\otimes AQ_1^{N_2}$ (see \cite{PP} for details). Since $E_1$ terms of spectral sequences coincide we conclude that $red$ is a quasi-isomorphism.
\end{proof}

The space $\left(AQ_{-N_1}^{0}\right)^*\otimes AQ_{1}^{N_2}$ will be used to compute cohomology $\Lambda^{\frac{\infty}{2}}(LH)$. We denote this space by $AQ^{\frac{\infty}{2}}=(AQ_{-N_1}^{N_2})^{\frac{\infty}{2}}$. The reader should keep in mind that $\left(AQ_{-N_1}^{0}\right)^*\otimes AQ_{1}^{N_2}$ is a module over $AQ_{-N_1}^{0}\otimes AQ_{1}^{N_2}$, and operators of multiplication on $\lambda^i[s]\otimes \lambda^i[t]$, $-N_1\leq s\leq 0, 1\leq t\leq N_2$ are defined. 

\begin{proposition}\label{P:smallcomplex}

The cohomology of $\Lambda^{\frac{\infty}{2}}(LH)$ coincide with  the cohomology of a two-step complex 
\begin{equation}\label{E:redquadric}
AQ^{\frac{\infty}{2}} \overset{d}\rightarrow AQ^{\frac{\infty}{2}}
\end{equation}
The map $d$ is defined by the formula 
\begin{equation}\label{E:d}
d(a)=\sum_{s+t=1, s\leq 0, t\geq 1} \lambda^i[s]\otimes \lambda^i[t]a , a\in AQ^{\frac{\infty}{2}} 
\end{equation}
The indices $s,t$ are in the range defined above.
\end{proposition}
\begin{proof}

The complex $ \Lambda^{\frac{\infty}{2}}(LH)_{red}$ as a linear space is isomorphic to 
\[AQ^{N_2}_{1}\otimes \Lambda(LH_{\leq 0})\otimes  \langle \omega \rangle\cong  AQ^{N_2}_{1}\otimes \Lambda(LH_{-N_1}^0)\otimes (\langle b[1]  \omega\rangle+ \langle \omega \rangle)\]
In the last isomorphism we used decomposition $LH_{\leq 0}=LH_{-N_1}^0+\langle b[1] \rangle$.
Introduce a notation for the graded space \[\langle b[1]  \omega\rangle+ \langle \omega\rangle=B^{\frac{\infty}{2}-1}+B^{\frac{\infty}{2}+0}\]
We denote temporally $A=AQ^{N_2}_{1}$. Recall that $A$ is graded .
A larger graded space $B=\bigoplus B^{\frac{\infty}{2}+i}$ is defined as a sum\[B^{\frac{\infty}{2}+i}=A_i\otimes B^{\frac{\infty}{2}+0}+A_{i+1}\otimes B^{\frac{\infty}{2}-1} \]
Filtration  of the graded algebra $A$ \[F^k=\bigoplus_{k\geq i} A_i\] defines a filtration  on  $\Lambda^{\frac{\infty}{2}}(LH)_{red}$. The corresponding  spectral sequence has $E_1$ term :
\[E^{i,j}_1=H_{-i}(LH_{-N_1}^0,\mathbb{C})\otimes B^{\frac{\infty}{2}+i}\rightarrow H^{\frac{\infty}{2}+i+j}(LH)_{red} \]
Homology $H_{i}(LH_{-N_1}^0,\mathbb{C})$ is dual to cohomology , which we know how to compute. We conclude  that  $H_{\bullet}(LH_{-N_1}^0,\mathbb{C})=\left(AQ^{0}_{-N_1}\right)^*$.
The $E_2$-term coincides with (\ref{E:redquadric}) with the following identifications: element $a$ in the domain of $d$ gets identified with  \[a b[1]\omega\] in some sub-quotient of   $\Lambda^{\frac{\infty}{2}}(LH)$. Likewise the image of $d$ is an element  \[\sum_{s+t=1, s\leq 0, t\geq 1} \lambda^i[s]\otimes \lambda^i[t]a  \omega\] The spectral sequence degenerates in $E_3$ term because of the two-step nilpotency of the algebra $LH$. 

\end{proof}
\section{Limits $N_1,N_2\rightarrow \infty$}\label{S:limits}
Our goal is to define limits of $\Lambda^{\frac{\infty}{2}}(LH_{-N_1}^{N_2})$ and study its properties.
 Only in these  limits Virasoro and affine Lie algebra actions emerge . Taking these limits is not straightforward and we shall discuss this presently.

We start with a remark that there are embeddings \[LH_{-N_1}^{N_2}\subset LH_{-N'_1}^{N'_2}, N_1\leq N'_1, N_2\leq N'_2\]
It is well known that inclusions define maps in cohomology and homology complexes:
\[\begin{split}
&C^*(LH_{-N'_1}^{N'_2})\rightarrow C^*(LH_{-N_1}^{N_2})\\
&C_*(LH_{-N_1}^{N_2})\rightarrow  C_*(LH_{-N'_1}^{N'_2})
\end{split}\]
Semi-infinite cohomology is not a functor  of  the Lie algebra because  it shares properties of both homology and cohomology (see \cite{Voronov}).
It has  a weaker property:
\begin{proposition}
There are morphisms of complexes of $C^*(LH_{-N'_1}^{N_2})$ and $C^*(LH_{-N_1}^{N'_2})$ modules:
\[\begin{split}
&\Lambda^{\frac{\infty}{2}}(LH_{-N'_1}^{N_2})\rightarrow \Lambda^{\frac{\infty}{2}}(LH_{-N_1}^{N_2})\\
&\Lambda^{\frac{\infty}{2}}(LH_{-N_1}^{N_2})\rightarrow \Lambda^{\frac{\infty}{2}}(LH_{-N_1}^{N'_2})\\
\end{split}\]
$N_1\leq N'_1, N_2\leq N'_2$
\end{proposition}
\begin{proof}
The first map follows from the isomorphism of differential graded $C(LH_{-N'_1}^{N_2})$-modules  \[\Lambda^{\frac{\infty}{2}}(LH_{-N_1}^{N_2})\cong C^*(LH_{-N_1}^{N_2})\underset{C^*(LH_{-N'_1}^{N_2})}{\otimes } \Lambda^{\frac{\infty}{2}}(LH_{-N'_1}^{N_2})\]

The second from
  \[\Lambda^{\frac{\infty}{2}}(LH_{-N_1}^{N'_2})\cong C_*(LH_{-N_1}^{N'_2})\underset{C^*(LH_{-N_1}^{N'_2})}{\otimes } \Lambda^{\frac{\infty}{2}}(LH_{-N_1}^{N_2})\]
\end{proof}

We define $\Lambda^{\frac{\infty}{2}}(LH_{-\infty}^{\infty})$ to be a double limit
\[\Lambda^{\frac{\infty}{2}}(LH_{-\infty}^{\infty})=\underset{\underset{N_1}\longleftarrow}\lim\   \underset{\underset{N_2}\longrightarrow}\lim \Lambda^{\frac{\infty}{2}}(LH_{-N_1}^{N_2})\]

The space $PQ^{\frac{\infty}{2}}$ it the double limit
 \[PQ^{\frac{\infty}{2}}=\underset{\underset{N_1}\longleftarrow}\lim\   \underset{\underset{N_2}\longrightarrow}\lim AQ_{-N_1}^0\otimes AQ_{1}^{N_2 *}\]

\begin{proposition}
The cohomology of $\Lambda^{\frac{\infty}{2}}(LH_{-\infty}^{\infty})$ coincide with cohomology of a two-term complex
\begin{equation}\label{E:dinfty}
PQ^{\frac{\infty}{2}}\overset{d}\rightarrow PQ^{\frac{\infty}{2}}
\end{equation}
with the differential as in (\ref{E:d}).
\end{proposition}
\begin{proof}
First of all the operator $d$ is well defined because $\lambda^i[t]$ acts trivially on $ AQ_{1}^{N_2 *}$  with $t>N_2$.

Secondly, The complexes $\Lambda^{\frac{\infty}{2}}(LH_{-\infty}^{\infty})$ and (\ref{E:dinfty}) decompose into a direct product of eigenspaces of the dilation operator. The eigenspaces are finite-dimensional complexes. The spaces of these complexes stabilize for large $N_1$ and $N_2$. The cohomology can be computed using appropriate eigenspace of  $\Lambda^{\frac{\infty}{2}}(LH_{-N_1}^{N_2})$ and by Proposition \ref{P:smallcomplex} coincide with cohomology of the corresponding eigenspaces of (\ref{E:redquadric}).  The later eigenspaces stabilize  for $N'_1\geq N_1,N'_2\geq N_2$.
\end{proof}

There is a natural isomorphism $AQ_{-N_1}^{N_2}\cong AQ_{-N_2}^{N_1}$ defined by the formula $\lambda(z)\rightarrow \lambda(1/z)$ (compare with \cite{AA}). It enables us to define  an isomorphism 
\begin{equation}\label{E:duality}
\begin{split}
&\left(AQ_{-N_1}^{0*}\otimes AQ_{1}^{N_2 *}\right)^* \cong 
AQ_{-N_1}^0\otimes AQ_{1}^{N_2\  *}\cong\\ 
&AQ_{-N_1-1}^{-1}\otimes AQ_{0}^{N_2-1\  *}\cong 
AQ^{0\ *}_{-N_2+1 }\otimes AQ^{N_1+1}_{1}
\end{split}
\end{equation}

We also get a limiting isomorphism of cohomological degree one
\[PQ^{\frac{\infty}{2}}\cong PQ^{\frac{\infty}{2}*}\] 
We denote the pairing between linear spaces of the complex (\ref{E:dinfty}) by $<\cdot,\cdot>$.
\begin{lemma}
The map $d$ (\ref{E:dinfty}) satisfies $<d(a),b>=<a,d(b)>$.
\end{lemma}
\begin{proof}
 The key moment is that a shift by one in indices  in (\ref{E:duality}) agrees with the rule $s+t=1$ in (\ref{E:d}).
\end{proof}

As a  corollary we get the  following Proposition.
\begin{proposition}
The cohomology of the complex $\Lambda^{\frac{\infty}{2}}(LH_{-\infty}^{\infty})$ has a non-degenerate  Poincar\'{e} duality pairing. 
\end{proposition}
\begin{remark}
It is still remains unclear how to compute generating functions of individual cohomology groups  of  (\ref {E:redquadric}) or (\ref{E:dinfty}).
\end{remark}

{\bf\Large  Appendix}

\appendix

\section{Exceptional cases}
\subsection{$n=4$}
The case of two-dimensional projective quadric ( $n=4$ ) requires  additional notations.
The Lie algebra of the symmetry group of the quadric is non-simple  $\so_4\cong \sl_2\times\sl_2$. Choose  generators $a=a_+,b=b_+$ of nilpotent subalgebras in left and right copies of $\sl_2$. Elements $a_-$, $b_-$ are  such that $(a_+,[a_+,a_-],a_-),(b_+,[b_+,b_-],b_-) $ form $\sl_2$-triples.\footnote{Elements $(e,h,f)$ define an $\sl_2$-triple if $[h,e]=2e,[h,f]=-2f,[e,f]=h$.}Define \[\hat{a}=a_-\otimes z\in \sl_2\otimes \mathbb{C}[z,z^{-1}],\] \[\hat{b}=b_-\otimes z\in \sl_2\otimes \mathbb{C}[z,z^{-1}].\]It is useful to identify four-dimensional fundamental representation of $\so_4$ with the tensor product $V=W_l\otimes W_r$, where $W_l,W_r$ are spinorial representations of left and right copies of $\sl_2$. We choose  bases $(\xi_1,\xi_2)$ and $(\eta_1,\eta_2)$ of $W_l$, $W_r$ that consist of weight vectors. In this bases \[a(\xi_1)=\xi_2, b(\eta_2)=\eta_1\]
Linear spaces $W_l$, $W_r$ are equipped $\so_4$-invariant skew-symmetric dot-products $<\cdot,\cdot>$, $[\cdot,\cdot]$ normalized by conditions $<\xi_1,\xi_2>=1$, $[\eta_1,\eta_2]=1$. We set \[f^2=\xi_1\otimes \eta_2, f^1=\xi_2\otimes \eta_2\]\[g_1=\xi_1\otimes \eta_1, g_2=\xi_2\otimes \eta_1\]
The invariant dot-product on $V=W_l\otimes W_r$ is equal to $<\cdot,\cdot>\otimes[\cdot,\cdot]$.
 The Hasse diagram $\hat{\sG}(2)$ is 
\[\begin{array}{ccccccccc}
\dots \rightarrow&\langle f^2[t] \rangle&\overset{a}\rightarrow&\langle f^1[t] \rangle&\overset{\hat{a}}\rightarrow&\langle f^2[t+1] \rangle&\overset{a}\rightarrow &\langle f^1[t+1] \rangle&\rightarrow \dots\\
&\downarrow b&&\downarrow b&&\downarrow b&&\downarrow b\\
\dots \rightarrow&\langle g_1[t] \rangle&\overset{a}\rightarrow&\langle g_2[t] \rangle&\overset{\hat{a}}\rightarrow&\langle g_1[t+1] \rangle&\overset{a}\rightarrow &\langle g_2[t+1] \rangle&\rightarrow \dots
\end{array}\]
Additional arrows that did not fit to the diagram are
 \[\langle g_2[t]\rangle \overset{\hat{b}}\rightarrow \langle f^1[t+1] \rangle\]
 \[\langle g_1[t] \rangle\overset{\hat{b}}\rightarrow \langle f^2[t+1] \rangle\]
 This diagram {\it does not} define a lattice because, e.g. supremum of $g_1[t],f^1[t]$ consists of two elements: $g_2[t]$ and $f^2[t+1]$.
 
The diagram $\sG_{-N_1}^{N_2}$ coincides with $[\langle f^n[-N_1]\rangle,\langle g_n[N_2]\rangle]$.  
\begin{equation}\label{E:relation2}
\begin{split}
&g_1[t]f^1[t]=-g_2[t]f^2[t]-g_2[t-1]f^2[t+1]-\left(g_1[t-1]f^1[t+1]+\sum_{s\geq 2}\sum_{i=1}^2 g_i[t+s]f^i[t-s]\right)\\
&g_2[t]f^2[t+1]=-g_{1}[t+1]f^{1}[t]-g_{2}[t]f^{2}[t+1] -\left(g_1[t]f^1[t+1]+\sum_{s\geq 2}\sum_{i=1}^2 g_i[t+s]f^i[t+1-s]\right)
\end{split}
\end{equation}
\begin{proposition}
Suppose $n=4$.
We refine the partial order on the poset $\hat{\sG}(4)$ to the total order as follows:\[\cdots < \langle f^2[t] \rangle < \langle g_1[t]\rangle < \langle f^1[t]\rangle< \langle g_2[t] \rangle < \langle f^2[t+1]\rangle < \cdots\]
Then the generators (\ref{E:relation2}) of the ideal of relations defines a Gr\"{o}bner basis of the ideal with respect to the degree-lexicographic order on monomials.
\end{proposition}
\begin{proof}
The proof is the same as in the general case. 
\end{proof}

As a corollary we get that $AQ_{-N_1}^{N_2}(4)$ is a Hodge algebra(\cite{DCEP}, \cite{Hibi}).
\begin{proposition}
The generating function $AQ_{-N_1}^{N_2}(4)(q,t)$ , $,N_1,N_2\geq 0$ is 
\[\frac{\prod_{l=-2N_1}^{2N_2}(1-q^lt^2)}{\prod_{l=-N_1}^{N_2}(1-q^lt)^{4}}\]
\end{proposition}
\begin{proof}
Exercise.
\end{proof}
\subsection{$n=3$}
The fundamental representation $V$ of $\so_3\cong \sl_2$ coincides with the adjoint .  We choose $R$ to be a commutator with $g$, $\hat{R}$ to be commutator with $f\otimes z$. The periodic affine diagram is given below.
\[\begin{split}
& \hat{\sG}(3):\\
& \cdots \langle h[t] \rangle\begin{array}{c}\overset{\hat{R}\quad \text{ }}\nearrow  \begin{array}{c}\langle f[t+1] \rangle\\ \text{ } \end{array} \overset{\text{ }\quad R}\searrow \\ \underset{R \quad  \text{ }}\searrow   \begin{array}{c} \text{ }\\\langle g[t] \rangle \end{array} \underset{\text{ }\quad \hat{R}}\nearrow \end{array}\langle h[t+1] \rangle \cdots t\in \mathbb{Z}
\end{split}
 \]
We construct a weakly partly ordered set using the standard  prescription with a proviso that $\langle h[t] \rangle$ are not reflexive vertices. The attentive reader will immediately notice that $ \hat{\sG}(3)$ is not a weak lattice because  \[\sup\{\langle h[t] \rangle,\langle h[t] \rangle\}=\{\langle g[t] \rangle,\langle f[t+1] \rangle\}\]
 
\begin{proposition}
Suppose $n=3$.
We define a  total order on the set of generators of $AQ_{-N_1}^{N_2}(3)$ (see equation \ref{E:algebraodd} for definition ) as  follows:\[\cdots  < \langle g[t-1]\rangle < \langle f[t]\rangle<  \langle h[t]\rangle < \langle g[t] \rangle < \langle f[t+1]\rangle < \cdots\]
Then the generators (\ref{E:relation3}) of the ideal of relations defines a Gr\"{o}bner basis of the ideal with respect to the degree-lexicographic order on monomials.
 \begin{equation}\label{E:relation3}
\begin{split}
&(h[t])^2=-g[t]f[t]-g[t-1]f[t+1] -\sum_{s\neq 0} h[t+s]h[t-s] -\sum_{s\neq 0,-1}g[t+s]f[t-s]\\
&g[t]f[t+1]=-h[t]h[t+1]-\sum_{s\neq 0} \left(h[t+s]h[t+1-s]+ g[t+s]f[t+1-s]\right)
\end{split}
\end{equation}
In the above formulas only  variables $f[t],g[s],h[k]$ are present, that satisfy $-N_1\leq t,s,k\leq N_2$. 
\end{proposition}
\begin{proof}
The proof is no different then the proof of the general case.
\end{proof}
\begin{proposition}
The generating function $AQ_{-N_1}^{N_2}(3)(q,t)$ , $N_1,N_2\geq 0$ is 
\[\frac{\prod_{l=-2N_1}^{2N_2}(1-q^lt^2)}{\prod_{l=-N_1}^{N_2}(1-q^lt)^{3}}\]
\end{proposition}
\begin{proof}
Exercise.
\end{proof}
\subsection{$n=2$}

\begin{remark}{\rm
Suppose $n=2$.
We define a  total order on the set of generators of $AQ_{-N_1}^{N_2}(2)$ as  follows:\[\cdots < \langle f[t-1] \rangle < \langle g[t-1]\rangle < \langle f[t]\rangle< \langle g[t] \rangle < \langle f[t+1]\rangle < \cdots\]
Then the generators (\ref{E:relation1}) of the ideal of relations {\it do not }defines a Gr\"{o}bner basis of the ideal with respect to the degree-lexicographic order on monomials.
\begin{equation}\label{E:relation1}
\begin{split}
&g[t]f[t]=-\sum_{s\neq 0} g[t+s]f[t-s]\\
&g[t]f[t+1]=-\sum_{s\neq 0} g[t+s]f[t+1-s]\
\end{split}
\end{equation}
In the above formulas only  variables $f[t],g[s]$ are present, that satisfy $-N_1\leq t,s\leq N_2$. 
}
\end{remark}
The reason is that $S$-polynomials computed from the pair $g[t]f[t]$ and $g[t]f[t+1]$ are nontrivial. Nontrivial  also are reduction $A(S(g[t]f[t],g[t]f[t+1]))$.
 
 Here are results of computation of generating functions $AQ_{0}^{N}(2)(1,t)$ that utilize  {\it Maple} package {\it Gr\"{o}bner}:
 \begin{equation}
 \begin{split}
& AQ_{0}^{0}(2)(1,t)=\frac{t+1}{1-t}\\
& AQ_{0}^{1}(2)(1,t)=\frac{t^4-2t^3+2t+1}{(1-t)^2}\\
&AQ_{0}^{2}(2)(1,t)=\frac{-t^7-3t^6+11t^5-5t^4-5t^3+t^2+3t+1}{(1-t)^3} \\
&AQ_{0}^{3}(2)(1,t)=\frac{t^{10}+4t^9+3t^8-48t^7+56t^6-14t^4-8t^3+3t^2+4t+1}{(1-t)^4}\\
 \end{split}
 \end{equation}
 Algebra $AQ_{0}^{3}(2)$ is not Koszul because some coefficients of $AQ_{0}^{3}(2)(1,-t)^{-1}$ are negative.
 
 We also would like to bring readers attention to the fact that numerators of the above rational functions {\it are not} palindromic polynomials. This contrasts with palindromic property of numerators of (\ref{E:genfunct}) corresponding to a quadric in $n\geq3$ dimensional space.
\section{Gr\"{o}bner bases}\label{S:Grobner}
A  powerful technique of commutative Gr\"{o}bner bases 
significantly simplifies computation of Poincar\'{e} series of graded algebras. We will  review  following \cite{Mora} this technique in the present  section. 

Let $\mathbb{N}$ be the set of non-negative integers, $\mathbf{T}$ be a commutative semigroup generated by $a_1,\dots,a_n$, whose elements are $t=a_1^{e_1}\cdots a_n^{e_n}, e_i\in \mathbb{N}$. We define $\deg[t]=e_1+\cdots +e_n$.
We choose the following total degree-lexicographic well  order on $\mathbf{T}$:
\[t_1=a_1^{e_1}\cdots a_n^{e_n}<a_1^{f_1}\cdots a_n^{f_n}=t_2 \]
if and only if
\[\deg(t_1)<\deg(t_2) \mbox{ or } \deg(t_1)=\deg(t_2) \mbox{ and there is } j: i<j,e_i=f_i, e_j<f_j\]

The order is compatible with the semi-group operation in a sense that for $s,t_1,t_2\in \mathbf{T}$,  $t_1<t_2$ implies $st_1<st_2$.

The semigroup ring $\mathbb{C}[ \mathbf{T}]$ is nothing else but a polynomial algebra $\mathbb{C}[a_1,\dots,a_n]$.
Each element $f\in \mathbb{C}[ \mathbf{T}]$ has a unique ordered representation 
\begin{equation}\label{E:set}
f=\sum_{i=1}^s c_it_i, c_i\in \mathbb{C}^*, t_i\in \mathbf{T}, t_1>t_2\cdots >t_s 
\end{equation}

With  every nonzero element $f\in \mathbb{C}[\mathbf{T}]$ we can associate $T(f)=t_1$  - the maximal term of $f$ and $lc(f)=c_1$ - the leading coefficient of $f$.

If $I\subset \mathbb{C}[\mathbf{T}]$ is an ideal , the set \[T(I)=\{T(f)\in \mathbf{T}|f\in I \}\]
is a semigroup ideal. 
Introduce a set  \[E(I)=\mathbf{T}\backslash T(I)\]

If the ideal is graded, the generating function $E(I)[t]$ coincides  with Poincar\'{e} series of $A=\mathbb{C}[\mathbf{T}]/I$.

The following theorem holds:
\begin{theorem}\label{T:basic1}
\begin{enumerate}
\item  \label{I:one} The algebra $\mathbb{C}[\mathbf{T}]$  is a direct sum of $I$ and the span of $E(I)$: \begin{equation}\label{I:decomp} \mathbb{C}[\mathbf{T}]=I\oplus \langle E(I)\rangle.\end{equation}
\item There is a $\mathbb{C}$-vector space isomorphism $\mathbb{C}[\mathbf{T}]/I$ and $\langle E(I)\rangle$
\item Let $can:\mathbb{C}[\mathbf{T}]\rightarrow \langle E(I)\rangle$ be projection on the second summand in  (\ref{I:decomp}). The element $can(f,I)$ is the canonical form of  $f$ with respect to ideal $I$ and the degree-lexicographic order.  For each $f\in \mathbb{C}[\mathbf{T}]$ there is $g=can(f,I)\in \langle E(I)\rangle$ such that $f-g\in I$.
\end{enumerate}
Moreover 
\begin{enumerate}
\item $can(f,I)=can(g,I)$ if and only if $f-g\in I$
\item $can(f,I)=0$ if and only if $f\in I$
\end{enumerate}
\end{theorem}
Suppose we know all about the semi-group ideal $T(I)$. Then the canonical form $can(f,I)$ can be effectively computed.  The iterative  procedure that does this computation can be described as follows. 
\begin{procedure}\label{A:reduction}
Given a polynomial $f$ we  traverse through  the set  $\{t_i\}$ (\ref{E:set}) starting from the greatest element, descend down the order and  looking for an element $t=t_i\in T(I)$. As \[t=\sum_{g_k<t} c_kg_k  \mathrm{ mod  }\  I,\] we can replace $t$ in the sum (\ref{E:set})  by a linear combination  of monomials   of  lower order, still keeping the image of $f$ in $\mathbb{C}[\mathbf{T}]/I$ unchanged. After collecting similar terms we repeat the procedure.  Eventually  this procedure terminates. 
\end{procedure}

The reader should note that  
the described   modifications simplify the structure of $f$, because some monomials become smaller relative  to the order. In this sense $can(f,I)$ is the smallest representative in  the class $f$ mod $I$. We shall refer to substitution \[t \overset{r}\rightarrow \sum_{g_k<t} c_kg_k \] from Procedure \ref{A:reduction} as to reduction $r_t$.

In case $\mathbf{T}=\mathbb{N}$ we have $\mathbb{C}[\mathbf{T}]=\mathbb{C}[t]$; $\mathbb{C}[t]$ is a principal ideal domain  , i.e. for any ideal $I$ there is $g\in \mathbb{C}[t]$ such that \[I=(g).\]   The above procedure becomes the Euclidian division  algorithm  of polynomials in one variable: \[f=sg+r=sg+can(f,I(g))\]

One can think about this procedure  as a multidimensional generalization  of the Euclidian   algorithm (see also \cite{CLS}).

Remembering that each semigroup ideal in an ordered semigroup has a unique irredundant  basis  we obtain:
\begin{proposition}\cite{Mora}
If  $I\subset \mathbb{C}[\mathbf{T}]$ is an ideal, there is  a unique set $Gr(I)\subset I$ such that:
\begin{enumerate}
\item $G(I)=\{T(g)|g\in G\}$ is an irreducible basis of $T(I)$.
\item $lc(g)=1$ for each $g\in Gr$
\item $g=T(g)-can(T(g),I)$ for each $g\in Gr$
\end{enumerate}
$Gr$ is called the reduced Gr\"{o}bner basis of $I$.
\end{proposition}
 \begin{example}
 Let $I$ be an ideal in $\mathbb{C}[a_1,\dots,a_n]$ generated by the quadric   \[a^2_1+\dots +a^2_n.\]
 Then the semigroup ideal $T(I)$ is generated by $a^2_n$.
 This enables us to  compute, using the previous proposition as an aid,  the Poincar\'{e} series of the ideal \[I[t]=E(I)[t]=\frac{t^2}{(1-t)^{n}}\]
 \end{example}
 
 
 This example explains that the problem of computing Poincar\'{e} of an ideal can be reduced to a  combinatorial problem  of  Poincar\'{e} series of semigroup ideals.
 
Let $\{f_i\}$ a basis of $I$,i.e.  a minimal set of generators of  $I$.  Then \[T(\{f_i\})\subset G(I)\]
It is not true however that two sets should coincide.
\begin{example}
Let $f_1=x^3-2xy$ and $f_2 = x^2y-2y^2 +x$ Then
\[x (x^2y -2y^2 + x) - y (x^3 - 2xy) =x^2\]
so that $x^2\in I(f_1,f_2)$. Thus, $x^2 \in  T(I)$. However $x^2$ is not divisible by $T(f_1)=x^3$ or $T(f_2)=x^2y$ 
So $T(I(f_1,f_2))$ is not generated by $T(f_1),T(f_2)$

\end{example}

The procedure of finding $cal(f,I)$ relies on reductions $r_t$. In general the structure of $r_t, t\in G(I)$ could be quite irregular. It simplifies 
if   
\begin{equation}\label{E:grobnersimple}
G(I)=T(\{f_i\}).
\end{equation}   We  normalize $f_i$ by the condition  that 
\begin{equation}\label{E:lc}
lc(f_i)=1.
\end{equation}  Let $q_i\in \mathbf{T}$ be equal to $T(f_i)$. 

Procedure  (\ref{A:reduction}) is modified as follows
\begin{procedure}\label{A:reductionq}
Given a polynomial $f$ we  traverse through  the set  $\{t_i\}$ (\ref{E:set}) starting from the greatest element, descend down the order and  look for an element $t=t_i\in T(I)$. Under assumption (\ref{E:grobnersimple}) this means that $t$ is divisible by one of $q_j=q$: \begin{equation}\label{E:div}t=sq. \end{equation} This lets us to apply reduction $r_{sq}$, which  finishes an iteration of the algorithm. The algorithm terminates after a finite number of steps.
\end{procedure}
 Equation (\ref{E:div}) tells us that monomial $t$ contains $q$ as a sub-monomial. We may try to  find $can(f,I)$ by  repeatedly lowering  the order of monomials in $f$ by applying reductions to $q_i$-divisible monomials.
 \begin{theorem}\cite{Mora}
 Under  condition (\ref{E:grobnersimple})  algorithm \ref{A:reductionq} applied to $f$ gives $con(f,I)$.
 \end{theorem}
 The nontrivial part of this theorem is that   $con(f,I)$ is unambiguously defined. The issue is that some monomial $t$ in $f$ satisfies \[t=qs=q's'\mbox{ with } q'\neq q'\]
 This means that we can apply to distinct maximal sequences of reduction $r_{q_s}\cdots r_{q_1} $ ,  $r_{q'_s}\cdots r_{q'_1} $ to the element $f$. The  theorem asserts  that \[r_{q_s}\cdots r_{q_1} f=r_{q'_s}\cdots r_{q'_1}f= con(f,I).\]
 
 The procedure described in (\ref{A:reductionq}) can be formally  used even when condition (\ref{E:grobnersimple}) is not satisfied. In this case  however  \[r_{q_s}\cdots r_{q_1} f\neq r_{q'_s}\cdots r_{q'_1}f\neq con(f,I).\]

 This can  be rectified if we work with the set $G(I)$. It is convenient to enumerate elements of  $G$ by positive integers.
\begin{procedure}\label{A:reductiong}
Given a polynomial $f$ we traverse through  the set  $\{t_i\}$ (\ref{E:set}) starting from the greatest element, descend down the order and  look for an element \begin{equation}\label{E:div}t=sq, q\in G(I),s\in \mathbf{T} \end{equation} If factorization  $t=sq_i=s'q_j$ is not unique we choose  $q$ with the greater index. This lets us to apply reduction $r_{sq}$, which  finishes an iteration of the algorithm. The next iteration is applied to $r_{sq}f$. The procedure terminates after a finite number of steps and converges to $A(f,G)$.
\end{procedure}
  \begin{theorem}\cite{Mora}
 The procedure \ref{A:reductiong} applied to $f$ gives $A(f,G)=con(f,I)$. In particular $A(f,G)$ does not depend on how we enumerated elements of  $G$.
 \end{theorem}
 
 This theorem explains importance of sets  $G(I)$ and $Gr(I)$. The next construction known under the name of Buchberger algorithm  effectively  computes  $Gr(I)$.
 
 We still work under assumption (\ref{E:lc}). With  the notation $l.c.m(a,b)$ for the least common multiple for $a,b\in \mathbf{T}$ we define
\[T(f,g)=l.c.m(f,g)\]
\begin{definition}
The polynomial \[S(f,g)=\frac{T(f,g)}{T(g)}g-\frac{T(f,g)}{T(f)}f\]   is called $S$-polynomial of $f,g$.
\end{definition}
Observe that $T(\frac{T(f,g)}{T(g)}g)=T(\frac{T(f,g)}{T(f)}f)=T(f,g)$, which means that in the difference some cancelations should occur.

 The method of computation of $A(f,G)$ can be trivially extended to any set $F\supset \{f_i\} $. In this case we can not claim that $A(f,F)=con(f,I)$. The function $A(f,F)$ will be used for computation of Gr\"{o}bner basis $Gr(I)$.
 \begin{algorithm}
 Computation of $Gr(I)$
 
 As  zero approximation $Gr_0(I)$ to $Gr(I)$ we take a set $\{f_i\}$.
 
 We compute $Gr_i(I)$ as follows. 
 Suppose $Gr_{i-1}(I)=\{f_i|i=1,\dots n_{i-1}\}$.
 We compute all possible $S$-polynomials $S(f,f')$, $f,f'\in Gr_{i-1}(I) $ and then its reductions $A(S(f,f'),Gr_{i-1}(I))$ . Let $f_i, i=n_{i-1}+1,\dots, n_i$ be the set of nonzero reductions. Then $Gr_{i}(I)=\{f_i|i=1,\dots n_{i}\}$. 
 \end{algorithm}
 \begin{proposition}
 After finitely many steps the above algorithm terminates $Gr_{i-1}(I)=Gr_{i}(I)=Gr(I)$. Its result is a Gr\"{o}bner basis of $I$. 
 \end{proposition}
\begin{lemma}\label{L:prime}\cite{Mora} If $T(f_i,f_j)=T(f_i)T(f_j)$ then  $A(S(f_i,f_j),\{f_i\})=0$. 
\end{lemma}
This lemma  could be a useful, because it eliminates  computation of $A(S(f_i,f_j),\{f_i\})$ for a significant  set of  $S$-polynomials.


\begin{thebibliography}{99}
\bibitem{AA}Y. Aisaka, E. Aldo Arroyo   {\it Hilbert space of curved $\beta\gamma$ systems on quadric cones}	JHEP 0808:052,2008, 	arXiv:0806.0586v1 [hep-th]
\bibitem{AABN}Y. Aisaka, E. Aldo Arroyo, N. Berkovits, N. Nekrasov   {\it Pure Spinor Partition Function and the Massive Superstring Spectrum}	arXiv:0806.0584v1 [hep-th]
\bibitem{Bergman} G. Bergman   {\it The diamond lemma for rings theory.} Advances in Math. 29, N. 2, p. 178–218,
1978.
\bibitem{Berkovits}N. Berkovits {\it Super-Poincar\'{e} Covariant Quantization of the Superstring}	JHEP 0004:018,2000, 	arXiv:hep-th/0001035v2
\bibitem{Bezrukavnikov}R. Bezrukavnikov {\it Koszul Property and Frobenius Splitting of Schubert Varieties } 	arXiv:alg-geom/9502021v1
\bibitem{Bogvad}R. B\"{o}gvad {\it Some homogeneous coordinate rings that are Koszul algebras}	arXiv:alg-geom/9501011v2
\bibitem{Braverman} A. Braverman {\it Spaces of quasi-maps into the flag varieties and their applications} Proc. ICM 2006, Madrid, Zurich: Eur. Math. Soc., 2006 
\bibitem{Chevalley} C. Chevalley. {\it Theory of Lie groups I.} Princeton University Press, Princeton, NJ, 1933.
\bibitem{CLS}D. Cox, J. Little, D. O'Shea {\it Ideals, varieties, and algorithms: an introduction to computational algebraic geometry and commutative algebra} Springer, 2007 
\bibitem{DCEP} C. De Concini , D. Eisenbud, C. Procesi  {\it Hodge Algebras.} 
Astérisque 91 (1982)
\bibitem{Hibi}T. Hibi {\it Every affine graded ring has a Hodge algebra structure.} Rend. Sem. Mat. Univ. Polytech. Torino 44 (1986), 278-286
\bibitem{Hodge} W.V.D. Hodge {\it Some enumerative results in the theory of forms}, Proc. Camb. Phil. Soc. 39 (1943),22–30.
\bibitem{Kempf}  G. R. Kempf  {\it Some wonderful rings in algebraic geometry } Journal of Algebra,  Vol 134, Issue 1,October 1990, Pages 222-224 
\bibitem{Mora} T. Mora {\it An introduction to commutative and non-commutative Gr\"{o}bner bases}
\bibitem{Nekrasov} N. Nekrasov {\it Lectures on curved beta-gamma systems, pure spinors, and anomalies} arXiv:hep-th/0511008v1
\bibitem{PP}A. Polishchuk ,  L.  Positselski {\it Quadratic algebras} Univ. Lecture Ser., vol. 37, Amer. Math. Soc., Providence, RI, 2005
\bibitem{Priddy} S. Priddy {\it Koszul resolutions.} Trans. Amer. Math. Soc. 152, p. 39–60, 1970.
\bibitem{Ravi}M.S. Ravi {\it Coordinate rings of G/P are Koszul}  Journal of Algebra Volume 177, Issue 2, 1995, Pages 367-371 
\bibitem{Ruffo}J. Ruffo{\it Quasimaps, straightening laws, and quantum cohomology for the Lagrangian Grassmannian} arXiv:0806.0834v1  [math.AG]
\bibitem{SottileSturmfels} F. Sottile, B. Sturmfels, {\it A sagbi basis for the quantum Grassmannian,} J. Pure and Appl. Alg. 158 (2001), 347–366.
\bibitem{Voronov}A. A. Voronov   {\it Semi-infinite homological algebra}	Inventiones Mathematicae, Vol 113, N. 1, 1993






\end{thebibliography}
\end{document}